\documentclass[oneside,english]{amsart}
\usepackage[T1]{fontenc}
\usepackage[latin9]{inputenc}
\usepackage{amsthm}
\usepackage{graphicx}
\usepackage{amssymb}
\usepackage[unicode=true, 
 bookmarks=true,bookmarksnumbered=false,bookmarksopen=false,
 breaklinks=false,pdfborder={0 0 1},backref=false,colorlinks=false]
 {hyperref}

\providecommand{\tabularnewline}{\\}

\numberwithin{equation}{section} %% Comment out for sequentially-numbered
\numberwithin{figure}{section} %% Comment out for sequentially-numbered
\theoremstyle{plain}
\newtheorem{thm}{Theorem}[section]
  \theoremstyle{definition}
  \newtheorem{defn}[thm]{Definition}
  \theoremstyle{plain}
  \newtheorem{prop}[thm]{Proposition}
 \theoremstyle{definition}
  \newtheorem{example}[thm]{Example}
  \theoremstyle{plain}
  \newtheorem{cor}[thm]{Corollary}
  \theoremstyle{remark}
  \newtheorem{rem}[thm]{Remark}
  \theoremstyle{plain}
  \newtheorem{lem}[thm]{Lemma}

\usepackage{epsfig}
\newcommand{\intr}{\mbox{\rm int}}
\newcommand{\epi}{\mbox{\rm epi}}
\newcommand{\Tepi}{T_{\scriptsize{\mbox{\rm epi}(f)}}}
\newcommand{\rint}{\mbox{\rm rint}}
\newcommand{\conv}{\mbox{\rm conv}}

\newcommand{\spanv}{\mbox{\rm span}}

\begin{document}

\title{First order dependence on uncertainty sets in robust optimization}

\author{C.H. Jeffrey Pang}

\curraddr{Department of Combinatorics and Optimization, University of Waterloo,
200 University Ave West, Waterloo, ON, Canada N2L 3G1.}

\email{\href{mailto:chj2pang@math.uwaterloo.ca}{chj2pang@math.uwaterloo.ca}}



\date{\today}
\begin{abstract}
We show that a first order problem can approximate solutions of a
robust optimization problem when the uncertainty set is scaled, and
explore further properties of this first order problem.
\end{abstract}
\maketitle
\tableofcontents{}

\section{Introduction}

Robust optimization is the methodology of handling optimization problems
with uncertain data. In practice, the presence of uncertainties in
optimization problems can make nominal solutions meaningless. Such
uncertainties can come from data uncertainty in measurement and estimation,
or from uncertainty in implementation. We refer to the recent text
\cite{BGN09} for more details.

Consider the linear program: \begin{eqnarray*}
 & \min_{x} & \bar{c}^{T}x+\bar{d}\\
 & \mbox{s.t.} & \bar{A}x\leq\bar{b}.\end{eqnarray*}
To account for the uncertainties in the data $(\bar{A},\bar{b})$,
one instead considers a point $x$ to be feasible if it satisfies
\[
Ax\leq b\mbox{ for all }(A,b)\in\mathcal{U}.\]
Here, $\mathcal{U}$ is a set containing the nominal data $(\bar{A},\bar{b})$.
We can consider the translation $\Delta\mathcal{U}=\mathcal{U}-(\bar{A},\bar{b})$
and ask: What is the behavior of optimal solutions to the robust optimization
problem if the set $\Delta\mathcal{U}$ were to be scaled by some
factor $\epsilon$? A large value of $\epsilon$ corresponds to a
more robust solution, and a small value of $\epsilon$ places more
importance in the objective function. Understanding the dependence
of $\epsilon$ allows one to find a balance between optimization and
robustness. The first order dependence on $\epsilon$ is addressed
in Corollary \ref{cor:lin-approx} for linear programs and Theorem
\ref{thm:TFO-approx} for nonlinear programs.

The outline of this paper is as follows. We introduce robust linear
programming in Section \ref{sec:Robust-lin}. Before we introduce
robust nonlinear programming in Section \ref{sec:Robust-nonlin},
we recall some topics in variational analysis (or nonsmooth analysis)
as presented in the texts \cite{RW98,Mor06,Cla83} in Section \ref{sec:var-analysis}.
In Section \ref{sec:Robust-nonlin}, we also define the tangential
problem, which will be important in Theorem \ref{thm:TFO-approx},
our main result. We present first order properties of the tangential
problem in Section \ref{sec:First-order}, and study the effects of
sums of uncertainty sets in the tangential problem in Section \ref{sec:uncertain-sets}.

\section{\label{sec:Robust-lin}Robust linear programming}

We keep our presentation compatible with \cite{BGN09}, and begin
with the definition of the robust counterpart of a linear program.
\begin{defn}
(Robust counterpart) For $\bar{A}\in\mathbb{R}^{m\times n}$, $\bar{b}\in\mathbb{R}^{m}$,
$\bar{c}\in\mathbb{R}^{n}$ and $\bar{d}\in\mathbb{R}$, where $m\geq n$,
consider the linear program with parameters $(\bar{A},\bar{b},\bar{c},\bar{d})$
\begin{eqnarray}
 &  & \min_{x}\bar{c}^{T}x+\bar{d}\label{eq:original-LP}\\
 &  & \mbox{s.t. }\bar{A}x\leq\bar{b}.\nonumber \end{eqnarray}
The \emph{robust counterpart }(written as RC) of the above linear
program is \begin{eqnarray*}
 &  & \min_{x}\left\{ \hat{c}(x)=\sup_{(A,b,c,d)\in\mathcal{U}}[c^{T}x+d]\mid Ax\leq b\mbox{ for all }(A,b,c,d)\in\mathcal{U}\right\} \\
 & = & \min_{x,t}\{t\mid t\geq c^{T}x+d\mbox{, }Ax\leq b\mbox{ for all }(A,b,c,d)\in\mathcal{U}\},\end{eqnarray*}
where $\mathcal{U}$ is an uncertainty set for the parameters $(A,b,c,d)$,
with $(\bar{A},\bar{b},\bar{c},\bar{d})\in\mathcal{U}$.
\end{defn}
In a typical linear program, the variable $d$ does not affect the
minimizer, but one has to take perturbations in $d$ into account
in a robust optimization problem. The second formulation in the RC
shows that we can rewrite the linear program so that $c$ stays constant
at $\bar{c}$ and $d=0$. This is the approach we will take for the
rest of this section, and we define $\mathcal{U}$ to be a set containing
elements of the form $(A,b)$, where $(A,b)$ are close enough to
$(\bar{A},\bar{b})$. For more details, we refer to \cite{BGN09}.

We define $\Delta A$, $\Delta b$ and the set $\Delta\mathcal{U}$
by the relations\begin{eqnarray*}
\Delta A & := & A-\bar{A},\\
\Delta b & := & b-\bar{b},\\
\mbox{ and }\Delta\mathcal{U} & := & \mathcal{U}-(\bar{A},\bar{b}).\end{eqnarray*}
The vector $x$ can be chosen so that it stays feasible under these
first order perturbations. We write $x=\bar{x}+\Delta x$. The RC
is therefore simplified to \begin{eqnarray}
 &  & \min_{\Delta x}c^{T}(\bar{x}+\Delta x)+d\nonumber \\
 &  & \mbox{s.t. }[\bar{A}_{i}+\Delta A_{i}](\bar{x}+\Delta x)\leq[\bar{b}_{i}+\Delta b_{i}]\mbox{ for all }(\Delta A_{i},\Delta b_{i})\in\Delta\mathcal{U}_{i},\label{eq:RC-final}\end{eqnarray}
where $\Delta\mathcal{U}_{i}$ is the uncertainty in the $i$th row. 

When $\Delta\mathcal{U}$ is a small set , we seek to use a first
order approximation to determine a robustly feasible $x$. Letting
$x=\bar{x}+\Delta x$, and removing the second order term $(\Delta A_{i})(\Delta x)$
in \eqref{eq:RC-final} gives\begin{eqnarray}
\bar{A}_{i}(\bar{x}+\Delta x)+(\Delta A_{i})\bar{x}-(\bar{b}_{i}+\Delta b_{i}) & \leq & 0\mbox{ for all }(\Delta A_{i},\Delta b_{i})\in\Delta\mathcal{U}_{i},\label{eq:robust-feasibility}\\
\mbox{or }[\bar{A}_{i}\bar{x}-\bar{b}]+\bar{A}_{i}(\Delta x)+[(\Delta A_{i})\bar{x}-\Delta b_{i}] & \leq & 0\mbox{ for all }(\Delta A_{i},\Delta b_{i})\in\Delta\mathcal{U}_{i}.\nonumber \end{eqnarray}
If $\bar{A}_{i}\bar{x}-\bar{b}_{i}<0$ and if $\Delta\mathcal{U}$
were small enough, this constraint will not be tight in the optimization
problem. With these in mind, we define the first order problem of
a linear program.
\begin{defn}
\label{def:linear-FO-robust}(First order problem) Let $\bar{x}$
be an optimal solution to \eqref{eq:original-LP}. The \emph{first
order problem} is the problem\begin{eqnarray}
 &  & \min_{\gamma}c^{T}\gamma\label{eq:FO-feasible}\\
 &  & \mbox{s.t. }\bar{A}_{i}\gamma+[(\Delta A_{i})\bar{x}-\Delta b]\leq0\nonumber \\
 &  & \mbox{ for all }(\Delta A_{i},\Delta b_{i})\in\Delta\mathcal{U}_{i}\mbox{ whenever }\bar{A}_{i}\bar{x}-\bar{b}_{i}=0.\nonumber \end{eqnarray}
The first order problem can also be written as \begin{eqnarray*}
 &  & \min_{\gamma}c^{T}\gamma\\
 &  & \mbox{s.t. }\bar{A}_{i}\gamma\leq-\max_{(\Delta A,\Delta b)\in\Delta\mathcal{U}}[(\Delta A_{i})\bar{x}-\Delta b_{i}]\\
 &  & \mbox{ whenever }\bar{A}_{i}\bar{x}-\bar{b}_{i}=0.\end{eqnarray*}

\end{defn}
 In the case where the optimal solution $\bar{x}$ is nondegenerate,
i.e., when $B=\{i\mid\bar{A}_{i}\bar{x}=\bar{b}_{i}\}$ is of size
$n$ and $\bar{A}_{B}$ is invertible, the optimal solution $\bar{\gamma}$
of the first order problem is just $\bar{\gamma}=\bar{A}_{B}^{-1}w$,
where $w$ is the vector\begin{equation}
w=\left(\begin{array}{c}
-\max_{(\Delta A,\Delta b)\in\Delta\mathcal{U}}[(\Delta A_{i_{1}})\bar{x}-\Delta b_{i_{1}}]\\
\vdots\\
-\max_{(\Delta A,\Delta b)\in\Delta\mathcal{U}}[(\Delta A_{i_{n}})\bar{x}-\Delta b_{i_{n}}]\end{array}\right),\label{eq:LP-easy}\end{equation}
where $i_{1},\dots i_{n}$ are the $n$ elements in $B$. When $\bar{x}$
is a degenerate solution, the first order problem is still easy to
solve. We illustrate with a particular example that the tangential
constraints are easily obtained for rectangular uncertainty sets.

\begin{example}
(Rectangular uncertainty) Suppose that the uncertainty set $\Delta\mathcal{U}$
is rectangular, that is\begin{eqnarray*}
\Delta\mathcal{U} & := & \big\{(\Delta A,\Delta b):|\Delta A_{j,k}|\leq\epsilon_{j,k}\mbox{, }|\Delta b_{j}|\leq\delta_{j}\\
 &  & \qquad\mbox{ for all }j\in\{1,\dots,m\}\mbox{, }k\in\{1,\dots,n\}\big\}.\end{eqnarray*}
Then for each $i\in B$, \[
\max\left\{ (\Delta A_{i})\bar{x}-\Delta b_{i}\mid(\Delta A_{i},\Delta b_{i})\in\Delta\mathcal{U}_{i}\right\} =\delta_{i}+\sum_{k=1}^{n}\epsilon_{i,k}|\bar{x}_{k}|.\]

\end{example}
In Theorem \ref{thm:TFO-approx}, we will discuss how an adapted
first order problem gives a first order approximation of the solution
to a robust optimization problem in a general setting of nonlinear
programs. For now, we shall present the corollary in the simpler setting
of linear programming.
\begin{cor}
\label{cor:lin-approx}(to Theorem \ref{thm:TFO-approx}) (First order
approximation in linear programming)  Consider the robust optimization
problem\begin{eqnarray}
 &  & \min_{x}c^{T}x+d\nonumber \\
 &  & \mbox{s.t. }(\bar{A}_{i}+\Delta A_{i})x\leq(\bar{b}_{i}+\Delta b_{i})\mbox{ for all }(\Delta A_{i},\Delta b_{i})\in\epsilon\Delta\mathcal{U}_{i}\mbox{ for all }i,\label{eq:small-e-RO}\end{eqnarray}
and the first order problem \begin{eqnarray}
 &  & \min_{\gamma}c^{T}\gamma\nonumber \\
 &  & \mbox{s.t. }\bar{A}_{i}\gamma\leq-\max_{(\Delta A_{i},\Delta b_{i})\in\Delta\mathcal{U}_{i}}[(\Delta A_{i})\bar{x}-\Delta b_{i}]\label{eq:FO_approx_lin}\\
 &  & \phantom{\mbox{s.t. }}\mbox{ for all }i\mbox{ s.t. }\bar{A}_{i}\bar{x}=\bar{b}_{i}.\nonumber \end{eqnarray}
Let $\bar{\Gamma}$ be the set of optimal solutions to \eqref{eq:FO_approx_lin}.
Suppose 
\begin{enumerate}
\item $\Delta\mathcal{U}_{i}$ are compact convex sets.
\item $\bar{\Gamma}$ is bounded.
\item There is some $\gamma^{\prime}$ such that $\bar{A}_{i}\gamma^{\prime}<0$
whenever $\bar{A}_{i}\bar{x}=\bar{b}_{i}$.
\item $\bar{x}$ is the unique minimizer of the nominal problem $\min\{c^{T}x\mid\bar{A}x\leq\bar{b}\}$.
\end{enumerate}
Then the set of cluster points of any sequence $\{\frac{1}{\epsilon}(\bar{x}_{\epsilon}-\bar{x})\}$,
where $\bar{x}_{\epsilon}$ is an optimal solution to \eqref{eq:small-e-RO}
and $\epsilon\to0$, is a subset of $\bar{\Gamma}$. The objective
value of \eqref{eq:small-e-RO} , say $\bar{v}_{\epsilon}$, has an
approximation $\bar{v}_{\epsilon}=\bar{v}+\epsilon\tilde{v}+o(\epsilon)$,
where $\tilde{v}$ is the objective value of \eqref{eq:FO_approx_lin}.

In particular, if $\bar{\Gamma}$ contains only one element, say $\bar{\gamma}$,
then $\lim_{\epsilon\to0}\frac{1}{\epsilon}(\bar{x}_{\epsilon}-\bar{x})=\bar{\gamma}$,
or $\bar{x}_{\epsilon}\in\bar{x}+\epsilon\bar{\gamma}+o(\epsilon)$.\end{cor}
\begin{proof}
The condition that $Q$ is Clarke regular at $\bar{A}\bar{x}-\bar{b}$
holds in this case because $\mathbb{R}_{-}^{m}$ is Clarke regular
everywhere. The condition that $\bar{x}$ is the unique minimizer
in (4) suffices because the domain is convex. The affine function
$x\mapsto c^{T}x+d$ is locally Lipschitz and subdifferentially regular
everywhere.
\end{proof}

\section{\label{sec:var-analysis}Preliminaries in variational analysis}

In this section, we recall the definitions of some nonsmooth objects
in variational analysis that will be necessary for the rest of the
paper. We recall the definition of normal cones and Clarke regularity.
\begin{defn}
\label{def:normal}(Normal cones and Clarke regularity) Let $C\subset\mathbb{R}^{n}$.
For a point $\bar{x}\in C$, a vector $v$ is \emph{normal to $C$
at $\bar{x}$ in the regular sense, }or a \emph{regular normal}, written
$v\in\hat{N}_{C}(\bar{x})$, if \[
v^{T}(x-\bar{x})\leq o(|x-\bar{x}|)\mbox{ for all }x\in C.\]
It is \emph{normal to $C$ in the general sense}, or simply a \emph{normal
vector}, written $v\in N_{C}(\bar{x})$, if there are sequences $x_{i}\to\bar{x}$
and $v_{i}\to v$ with $v_{i}\in\hat{N}_{C}(x_{i})$. The set $C$
is \emph{Clarke regular }at $\bar{x}$ if $N_{C}(\bar{x})=\hat{N}_{C}(\bar{x})$.
\end{defn}
We refer the reader to \cite[Corollary 6.29]{RW98} for equivalent
definitions of Clarke regularity. The sets we will encounter in this
paper are all Clarke regular, so this does not cause difficulties.

We recall the definition of the tangent cone, which will be important
in our main result.
\begin{defn}
\label{def:tangent}(Tangent cones) The \emph{tangent cone }of a set
$C\subset\mathbb{R}^{m}$ at some $\bar{x}\in C$ is defined by\[
T_{C}(\bar{x}):=\left\{ w\mid\frac{x_{i}-\bar{x}}{t_{i}}\to w\mbox{ for some }x_{i}\in C\mbox{, }t_{i}\searrow0\mbox{ and }x_{i}\to\bar{x}\right\} .\]

\end{defn}
Next, we recall sublinearity and equivalent definitions of subdifferential
regularity that will also be useful for our main result. We take the
definitions of subdifferential regularity from \cite[Definition 7.25, Exercise 9.15, Corollary 8.19]{RW98}.
\begin{defn}
(positive homogeneity and sublinearity) A function $h:\mathbb{R}^{n}\to\mathbb{R}$
is \emph{positively homogeneous} if $h(\lambda x)=\lambda h(x)$ for
all $x$ and $\lambda>0$. It is \emph{sublinear} if in addition \[
h(x+x^{\prime})\leq h(x)+h(x^{\prime})\mbox{ for all }x\mbox{ and }x^{\prime}.\]
It is clear that sublinear functions are convex.
\begin{defn}
\label{def:equiv-subdif-reg}(Subdifferential regularity) Let $f:\mathbb{R}^{n}\to\mathbb{R}$
be locally Lipschitz at $\bar{x}$. 

(a) We say that the function $f$ is \emph{(subdifferentially) regular
at $\bar{x}$ }if the epigraph $\epi f:=\{(x,t)\mid t\geq f(x)\}$
is Clarke regular at $(x,f(x))$ as a subset of $\mathbb{R}^{n}\times\mathbb{R}$.

(b) Define the \emph{subderivative} $df(\bar{x}):\mathbb{R}^{n}\to\mathbb{R}$
by \begin{equation}
df(\bar{x})(w):=\liminf_{\tau\searrow0}\frac{f(\bar{x}+\tau w)-f(\bar{x})}{\tau}.\label{eq:subderivative}\end{equation}
and the \emph{regular subderivative }$\hat{d}f(\bar{x}):\mathbb{R}^{n}\to\mathbb{R}$
by\[
\hat{d}f(\bar{x})(w):=\limsup_{\scriptsize{\begin{array}{c}
\tau\searrow0\\
x\to\bar{x}\end{array}}}\frac{f(x+\tau w)-f(x)}{\tau}.\]
In general, the regular subderivative is sublinear. The function $f$
is (subdifferentially) regular at $\bar{x}$ if and only if $df(\bar{x})=\hat{d}f(\bar{x})$.
Under subdifferential regularity, it is clear that the liminf in \eqref{eq:subderivative}
can be taken to be a full limit. Also, $\Tepi(\bar{x},f(\bar{x}))=\epi(df(\bar{x}))$.
\end{defn}
\end{defn}
Since the tangent cone will play a major role in our main result,
we now recall some calculus rules for tangent cones, highlighting
a constraint qualification condition similar to that of condition
\eqref{enu:CQ} in Theorem \ref{thm:TFO-approx}. The rest of this
section will not be essential to the development of the paper, so
one may skip to the next section in a first reading. We now recall
a formula for tangent cones under intersections.
\begin{prop}
\label{pro:tangent-NMFCQ}(Tangent cones to intersections) Let $C=C_{1}\cap\cdots\cap C_{m}$
for closed sets $C_{i}\subset\mathbb{R}^{n}$, and let $\bar{x}\in C$.
Suppose $\bar{x}$ is Clarke regular at $C_{j}$ for all $j$. Assume
either\begin{eqnarray}
 &  & \sum_{j=1}^{m}\lambda_{j}v_{j}=0\mbox{, }v_{j}\in N_{C_{j}}(\bar{x})\mbox{ and }\lambda_{j}\geq0\mbox{ for all }j\in\{1,\dots,m\}\label{eq:NMFCQ1}\\
 &  & \mbox{ implies }\lambda_{j}=0\mbox{ for all }j\in J^{\prime},\nonumber \end{eqnarray}
or equivalently:
\begin{enumerate}
\item [(a)]there are no vectors $\{y_{j}\}_{j=1}^{m}$ such that $y_{j}\perp T_{C_{j}}(\bar{x})$
and $y_{1}+\cdots+y_{m}=0$ other than $y_{j}=0$ for all $j\in\{1,\dots,m\}$,
and there is a vector $w$ such that $w\in\mathbb{R}^{n}\backslash\{0\}$
such that $w\in\rint(T_{C_{j}}(\bar{x}))$ for all $j\in\{1,\dots,m\}$. 
\end{enumerate}
Then one has\[
T_{C}(\bar{x})=T_{C_{1}}(\bar{x})\cap\cdots\cap T_{C_{m}}(\bar{x}),\]
and $C$ is Clarke regular at $\bar{x}$.\end{prop}
\begin{proof}
Other than the equivalence of \eqref{eq:NMFCQ1} and (a), this result
is stated in a more general case in \cite[Theorem 6.42]{RW98}. This
result is obtained by consider the set $D:=C_{1}\times\cdots\times C_{2}\subset(\mathbb{R}^{n})^{m}$
and the mapping $F:x\mapsto(x,\dots,x)\in(\mathbb{R}^{n})^{m}$ with
$X=\mathbb{R}^{n}$ and applying \cite[Theorems 6.31 and 6.41]{RW98}.
The constraint qualification condition required is \eqref{eq:NMFCQ1}.
By \cite[Exercise 6.39(b)]{RW98}, \eqref{eq:NMFCQ1} is equivalent
to the existence of a $w^{\prime}$ such that $F(w^{\prime})\in\rint(T_{D}(\bar{x},\dots,\bar{x}))$
and having \[
y\perp T_{D}(\bar{x},\dots,\bar{x})\mbox{ and }F^{*}(y)=0\mbox{ implies }y=0.\]
These conditions are equivalent to that in (a).
\end{proof}
We recall the Mangasarian-Fromovitz constraint qualification.
\begin{defn}
\label{def:MFCQ}(Mangasarian-Fromovitz constraint qualification)
For $\mathcal{C}^{1}$ functions $f_{j}:\mathbb{R}^{n}\to\mathbb{R}$
and $j\in\{1,\dots,m\}$, let \[
Q:=\big\{x\in\mathbb{R}^{n}\mid f_{j}(x)\leq0\mbox{ for all }j\in\{1,\dots,m\}\big\}.\]
For $\bar{x}\in C$, let $J^{\prime}:=\{j\mid f_{j}(\bar{x})=0\}$.
The \emph{Mangasarian-Fromovitz constraint qualification} (MFCQ) is
satisfied at $\bar{x}$ if there is a vector $w\in\mathbb{R}^{n}$
such that \[
\nabla f_{j}(\bar{x})^{T}w<0\mbox{ for all }j\in J^{\prime}.\]
Another equivalent definition of the MFCQ is the following {}``positive
linear independence'' condition \[
\sum_{j\in J^{\prime}}\lambda_{j}\nabla f_{j}(\bar{x})=0\mbox{ and }\lambda_{j}\geq0\mbox{ for all }j\in J^{\prime}\mbox{ implies }\lambda_{j}=0\mbox{ for all }j\in J^{\prime}.\]

\end{defn}
The classical definition of the MFCQ also takes into account equality
constraints in the set $Q$, which we omit since they are not of immediate
interest. 

To handle sets defined by nonsmooth constraints, we need to recall
the subdifferential.
\begin{defn}
(Subdifferentials)\label{def:subdifferential} Consider a function
$f:\mathbb{R}^{n}\rightarrow\mathbb{R}$ such that $f$ is locally
Lipschitz at $\bar{x}$. For a vector $v\in\mathbb{R}^{n}$, one says
that

(a) $v$ is a \emph{regular subgradient} (also known as a \emph{Fr\'{e}chet}
\emph{subgradient}) of\emph{ $f$ }at $\bar{x}$, written $v\in\hat{\partial}f(\bar{x})$,
if \[
f(x)\geq f(\bar{x})+\left\langle v,x-\bar{x}\right\rangle +o(\left|x-\bar{x}\right|);\]

(b) $v$ is a \emph{(general) subgradient} of $f$ at $\bar{x}$,
written $v\in\partial f(\bar{x})$, if there are sequences $x_{i}\rightarrow\bar{x}$
and $v_{i}\to v$ such that $f(x_{i})\rightarrow f(\bar{x})$ and
$v_{i}\in\hat{\partial}f(x_{i})$.

(c) The set $\hat{\partial}f(\bar{x})$\emph{ }is the \emph{regular
subdifferential,} and the set $\partial f(\bar{x})$ is the \emph{(general)
subdifferential}.

(d) The function $f$ is (subdifferentially) regular at $\bar{x}$\emph{
}if and only if $\partial f(\bar{x})=\hat{\partial}f(\bar{x})$.
\end{defn}
This characterization of subdifferentially regular functions is slightly
different from the earlier definitions, but is equivalent in the case
of locally Lipschitz functions in view of \cite[Corollary 8.11, Theorem 9.13 and Theorem 8.6]{RW98}.
We shall only be concerned with subdifferentially regular functions
throughout this paper, so there is no need to distinguish between
$\partial f(\bar{x})$ and $\hat{\partial}f(\bar{x})$. We conclude
with results on the intersections of tangent cones described by constraints.
\begin{prop}
\label{pro:tangent-examples}(Tangent cone under constraints) Suppose
$C=\{x\mid f_{j}(x)\leq0,j\in J\}$, and $J$ is a finite set. At
the point $\bar{x}\in C$, let $J^{\prime}\subset J$ be the set of
all $j$'s such that $f_{j}(\bar{x})=0$. If $f_{j}$ are continuous
at $\bar{x}$ for all $j\in J$, $f_{j}$ are continuously differentiable
at $\bar{x}$ for all $j\in J^{\prime}$ and the MFCQ is satisfied
at $\bar{x}\in C$, then\[
T_{C}(\bar{x})=\{z\mid\nabla f_{j}(\bar{x})^{T}z\leq0\mbox{ for all }j\in J^{\prime}\}.\]
In the nonsmooth case, if $f_{j}$ were locally Lipschitz and subdifferentially
regular at $\bar{x}$ for all $j\in J^{\prime}$ and \begin{eqnarray}
 &  & \sum_{j\in J^{\prime}}\lambda_{j}v_{j}=0\mbox{, }v_{j}\in\partial f_{j}(\bar{x})\mbox{ and }\lambda_{j}\geq0\mbox{ for all }j\in J^{\prime}\label{eq:nonsmooth-MFCQ}\\
 &  & \mbox{ implies }\lambda_{j}=0\mbox{ for all }j\in J^{\prime},\nonumber \end{eqnarray}
then $\bar{x}$ is Clarke regular at $C$, and \begin{eqnarray}
T_{C}(\bar{x}) & = & \bigcap_{j\in J^{\prime}}\{z\mid v^{T}z\leq0\mbox{ for all }v\in\partial f_{j}(\bar{x})\}.\label{eq:Tangent_Qi_formula}\\
 & = & \{z\mid v^{T}z\leq0\mbox{ for all }v\in\bigcup_{j\in J^{\prime}}\partial f_{j}(\bar{x})\}.\nonumber \end{eqnarray}
\end{prop}
\begin{proof}
We prove the general nonsmooth case for this theorem, which implies
the smooth case. There is a neighborhood $U$ of $\bar{x}$ such that
$C\cap U=[\cap_{j\in J^{\prime}}C_{j}]\cap U$, where $C_{j}$ is
defined by $C_{j}=\{x\mid f_{j}(x)\leq0\}$. Furthermore, $0\notin\partial f_{j}(\bar{x})$
for all $j\in J^{\prime}$. By \cite[Theorem 10.3]{RW98} (normal
cones to level sets) and \cite[Corollary 6.29(d)]{RW98} (tangent-normal
relations in regular sets), $C_{j}$ is Clarke regular at $\bar{x}$,
and the tangent cones $T_{C_{j}}(\bar{x})$ and normal cones $N_{C_{j}}(\bar{x})$
are given by \begin{eqnarray*}
N_{C_{j}}(\bar{x}) & = & \{\lambda v\mid\lambda\geq0\mbox{, and }v\in\partial f_{j}(\bar{x})\},\\
\mbox{ and }T_{C_{j}}(\bar{x}) & = & \{w\mid w^{T}v\leq0\mbox{ for all }v\in\partial f_{j}(\bar{x})\}.\end{eqnarray*}
Therefore, condition \eqref{eq:nonsmooth-MFCQ} becomes\[
\sum_{j\in J^{\prime}}v_{j}=0\mbox{, }v_{j}\in N_{C_{j}}(\bar{x})\mbox{ implies }v_{j}=0\mbox{ for all }j\in J^{\prime}.\]
By Proposition \ref{pro:tangent-NMFCQ}, the tangent cone $T_{C}(\bar{x})$
is \begin{eqnarray*}
T_{C}(\bar{x}) & = & \bigcap_{i\in J^{\prime}}T_{C_{j}}(\bar{x}),\end{eqnarray*}
which gives the formula for the tangent cone in the statement. 
\end{proof}
It is well known that for the sets \begin{eqnarray*}
C_{1} & := & \big\{(x_{1},x_{2})\in\mathbb{R}^{2}:x_{2}\geq x_{1}^{2}\big\}\\
C_{2} & := & \big\{(x_{1},x_{2})\in\mathbb{R}^{2}:x_{2}\leq-x_{1}^{2}\big\},\end{eqnarray*}
we have $T_{C_{1}\cap C_{2}}(0)\subsetneq T_{C_{1}}(0)\cap T_{C_{2}}(0)$
but the MFCQ is not satisfied.

The constraint qualification condition \eqref{eq:nonsmooth-MFCQ}
can be checked by another equivalent condition when $T_{C_{j}}(\bar{x})$
have nonempty interior.
\begin{prop}
\label{pro:CQ}(Constraint qualification) Assume the conditions of
Proposition \ref{pro:tangent-examples}. If $T_{C_{j}}(\bar{x})$
have nonempty interior and $0\notin\partial f_{j}(\bar{x})$ for all
$j\in J^{\prime}$, then the condition \eqref{eq:nonsmooth-MFCQ}
is equivalent to the existence of a vector $w$ such that $w\in\intr(T_{C_{j}}(\bar{x}))$
(or equivalently $w^{T}v<0$ for all $v\in\partial f_{j}(\bar{x})$)
for all $j\in J^{\prime}$. \end{prop}
\begin{proof}
Recall that by \cite[Theorem 10.3]{RW98}, if $0\notin\partial f(\bar{x})$,
then the tangent cone $T_{C_{j}}(\bar{x})$ is equal to $\{z\mid z^{T}v\leq0\mbox{ for all }v\in\partial f_{j}(\bar{x})\}$,
and the interior $\intr(T_{C_{j}}(\bar{x}))$ is $\{z\mid z^{T}v<0\mbox{ for all }v\in\partial f_{j}(\bar{x})\}$,
which gives the equivalence on the conditions on $w$. Next, the equivalence
of \eqref{eq:nonsmooth-MFCQ} and the condition in this result follow
from Proposition \ref{pro:tangent-NMFCQ}.
\end{proof}

\section{\label{sec:Robust-nonlin}Robust nonlinear programming}

We look at nonlinear programs of the form\begin{equation}
\min_{x}\{c^{T}x+d\mid\bar{A}x-\bar{b}\in Q\},\label{eq:RO_1}\end{equation}
where $Q\subset\mathbb{R}^{k}$ is a closed set. Specifically, we
consider problems of the form\begin{equation}
\min_{x}\{c^{T}x+d\mid\bar{A}_{i}x-\bar{b}_{i}\in Q_{i},\,1\leq i\leq m\},\label{eq:RO_2}\end{equation}
where $Q_{i}\subset\mathbb{R}^{k_{i}}$ are nonempty closed sets,
$\bar{A}_{i}\in\mathbb{R}^{k_{i}\times n}$, and $\bar{b}_{i}\in\mathbb{R}^{k_{i}}$.
We may write $\bar{A}$ as a concatenation of the matrices $\bar{A}_{i}$
and $\bar{b}$ as a concatenation of the vectors $\bar{b}_{i}$, and
this would make \eqref{eq:RO_1} equivalent to \eqref{eq:RO_2} for
$Q=Q_{1}\times\cdots\times Q_{m}$ and $k=k_{1}+\cdots+k_{m}$. One
case of interest is the set $Q_{i}=\{y\mid f_{i,j}(y)\leq0\mbox{ for all }j\in J\}$
for some $f_{i,j}:\mathbb{R}^{k_{i}}\to\mathbb{R}$ and the set $J$
is finite. Another case of interest is \emph{conic programs}, which
arise when all $Q_{i}$'s are closed convex pointed cones with nonempty
interior.

We now recall the definition of robust feasibility from \cite{BGN09}.
\begin{defn}
\label{def:Rob-feas}(Robust feasibility) Let an uncertain problem
be given and $\Delta\mathcal{U}=\Delta\mathcal{U}_{1}\times\cdots\times\Delta\mathcal{U}_{m}$
be a perturbation set. A candidate solution $x\in\mathbb{R}^{n}$
is \emph{robustly feasible} if it remains feasible for all realizations
of the perturbation vector from the perturbation set, that is \begin{equation}
[\bar{A}_{i}+\Delta A_{i}]x-[\bar{b}_{i}+\Delta b_{i}]\in Q_{i}\,\forall(i,1\leq i\leq m,(\Delta A_{i},\Delta b_{i})\in\Delta\mathcal{U}_{i}),\label{eq:conic-feasible}\end{equation}
where $\Delta\mathcal{U}_{i}\subset(\mathbb{R}^{k_{i}\times n}\times\mathbb{R}^{k_{i}})$
is the uncertainty set in $(\bar{A}_{i},\bar{b}_{i})$.\end{defn}
\begin{rem}
(Decomposing uncertainty sets) In the case where $\Delta\mathcal{U}$
is not a direct product of uncertainty sets, the uncertainty sets
$\Delta\mathcal{U}_{i}$ can be defined as \[
\Delta\mathcal{U}_{i}:=\{(\Delta A_{i},\Delta b_{i}):(\Delta A_{i},\Delta b_{i})=\Pi_{i}(\Delta A,\Delta b)\mbox{ for some }(\Delta A,\Delta b)\in\Delta\mathcal{U}\},\]
where $\Pi_{i}$ is the relevant projection from $\mathbb{R}^{k\times n}\times\mathbb{R}^{k}$
to $\mathbb{R}^{k_{i}\times n}\times\mathbb{R}^{k_{i}}$. It is clear
that \eqref{eq:conic-feasible} is equivalent to \[
[\bar{A}+\Delta A]x-[\bar{b}+\Delta b]\in Q\,\forall(\Delta A,\Delta b)\in\Delta\mathcal{U}.\]

\end{rem}
The definition for robust nonlinear programs encompasses nonlinear
objective functions.
\begin{example}
\label{exa:nonlin-obj}(Nonlinear objective) Consider the robust optimization
problem\begin{eqnarray*}
 &  & \min_{x}f(x)\\
 &  & \mbox{s.t. }[\bar{A}_{i}+\Delta A_{i}]x-[\bar{b}_{i}+\Delta b_{i}]\in Q_{i}\,\forall(i,1\leq i\leq m,(\Delta A_{i},\Delta b_{i})\in\Delta\mathcal{U}_{i}).\end{eqnarray*}
We can rewrite this robust problem as \begin{eqnarray*}
 &  & \min_{x,t}t\\
 &  & \mbox{s.t. }[\bar{A}_{i}+\Delta A_{i}]x-[\bar{b}_{i}+\Delta b_{i}]\in Q_{i}\,\forall(i,1\leq i\leq m,(\Delta A_{i},\Delta b_{i})\in\Delta\mathcal{U}_{i})\\
 &  & \mbox{and }f(x)\leq t.\end{eqnarray*}
The function $f$ is convex if and only if the epigraph $\epi(f)=\{(x,t)\mid f(x)\leq t\}$
is convex. Similarly, for a function $f$ locally Lipschitz at $\bar{x}$,
the function $f$ is subdifferentially regular at $\bar{x}$ if and
only if $\epi(f)$ is Clarke regular at $\bar{x}$. To prove our results
for nonlinear functions, we can prove the result for linear objective
functions and then appeal to the second formulation to obtain the
result we need. 
\end{example}
The formula in the robust optimization constraint can be rewritten
as \begin{eqnarray*}
 &  & [\bar{A}_{i}+\Delta A_{i}](\bar{x}+\Delta x)-[\bar{b}_{i}+\Delta b_{i}]\in Q_{i}\\
 & \iff & [\bar{A}_{i}\bar{x}-\bar{b}_{i}]+\bar{A}_{i}(\Delta x)+[(\Delta A_{i})\bar{x}-\Delta b_{i}]+(\Delta A_{i})(\Delta x)\in Q_{i}.\end{eqnarray*}
As in linear programming, we eliminate the second order term $(\Delta A)(\Delta x)$
to obtain a first order approximation. For nonlinear programs, we
also need to approximate the set $Q_{i}$ at $\bar{A}_{i}\bar{x}-\bar{b}_{i}$
by the tangential approximation $T_{Q_{i}}(\bar{A}_{i}\bar{x}-\bar{b}_{i})+[\bar{A}_{i}\bar{x}-\bar{b}_{i}]$
at $\bar{A}_{i}\bar{x}-\bar{b}_{i}$. This gives our definition of
the tangential problem.
\begin{defn}
\label{def:tangent-pblm}(Tangential problem) Let $\bar{x}$ be an
optimal solution to a nonlinear programming problem with parameters
$(\bar{A},\bar{b})$ so that $Q_{i}$ is Clarke regular at $\bar{A}_{i}\bar{x}-\bar{b}_{i}$
for all $i$. The \emph{tangential problem }to the robust optimization
problem obtained with constraints as explained in Definition \ref{def:Rob-feas}
is\begin{eqnarray*}
 &  & \min_{\gamma}c^{T}\gamma\\
 &  & \mbox{s.t. }[\bar{A}_{i}\bar{x}-\bar{b}_{i}]+\bar{A}_{i}\gamma+[(\Delta A_{i})\bar{x}-\Delta b_{i}]\in T_{Q_{i}}(\bar{A}_{i}\bar{x}-\bar{b}_{i})+[\bar{A}_{i}\bar{x}-\bar{b}_{i}]\\
 &  & \phantom{\mbox{s.t. }}\mbox{ for all }(i,1\leq i\leq m,(\Delta A_{i},\Delta b_{i})\in\Delta\mathcal{U}_{i}),\end{eqnarray*}
or equivalently\begin{eqnarray}
 &  & \min_{\gamma}c^{T}\gamma\label{eq:tangential-FO}\\
 &  & \mbox{s.t. }\bar{A}_{i}\gamma+[(\Delta A_{i})\bar{x}-\Delta b_{i}]\in T_{Q_{i}}(\bar{A}_{i}\bar{x}-\bar{b}_{i})\nonumber \\
 &  & \phantom{\mbox{s.t. }}\mbox{ for all }(i,1\leq i\leq m,(\Delta A_{i},\Delta b_{i})\in\Delta\mathcal{U}_{i}),\nonumber \end{eqnarray}
which is also equivalent to \begin{eqnarray*}
 &  & \min_{\gamma}c^{T}\gamma\\
 &  & \mbox{s.t. }\bar{A}_{i}\gamma+L_{i}(\Delta\mathcal{U}_{i})\in T_{Q_{i}}(\bar{A}_{i}\bar{x}-\bar{b}_{i})\\
 &  & \phantom{\mbox{s.t. }}\mbox{ for all }(i,1\leq i\leq m),\end{eqnarray*}
where $L_{i}:\mathbb{R}^{m_{i}\times n}\times\mathbb{R}^{m_{i}}\to\mathbb{R}^{m_{i}}$
is defined by $L_{i}(\Delta A_{i},\Delta b_{i})=(\Delta A_{i})\bar{x}-\Delta b_{i}$.
We call the corresponding constraints to the tangential problem the
\emph{tangential constraints}.\end{defn}
\begin{rem}
(Clarke regularity assumption) The assumption that each $Q_{i}$ is
Clarke regular at $\bar{A_{i}}\bar{x}-\bar{b}_{i}$ in Definition
\ref{def:tangent-pblm} comes about because the set $Q_{1}\times\cdots\times Q_{m}$
is Clarke regular at $\bar{A}\bar{x}-\bar{b}$ if and only if $Q_{i}$
is Clarke regular at $\bar{A}_{i}\bar{x}-\bar{b}_{i}$ for all $i$,
and in this case,\[
T_{Q_{1}\times\cdots\times Q_{m}}(\bar{A}\bar{x}-\bar{b})=T_{Q_{1}}(\bar{A}_{1}\bar{x}-\bar{b}_{1})\times\cdots\times T_{Q_{m}}(\bar{A}_{m}\bar{x}-\bar{b}_{m}).\]
(see \cite[Proposition 6.41]{RW98}.) This property makes the tangential
problem independent of how we decompose the set $Q$ as a direct product
of sets.
\end{rem}
We give some examples of tangential constraints.
\begin{example}
\label{exa:TFO}(Examples of tangential constraints) (a) When $\bar{A}_{i}\bar{x}-\bar{b}_{i}=0$
and $Q_{i}$ is a closed convex cone, then $T_{Q_{i}}(0)=Q_{i}$.
In this case, the corresponding tangential constraint is obtained
by just removing the second order term $(\Delta A_{i})(\Delta x)$.

(b) When $\bar{A}_{i}\bar{x}-\bar{b}_{i}\in\intr(Q_{i})$, then $T_{Q_{i}}(\bar{A}_{i}\bar{x}-\bar{b}_{i})=\mathbb{R}^{k_{i}}$
and the corresponding tangential constraint vanishes. 
\end{example}
In view of Example \ref{exa:TFO}, we see that for linear programming,
the tangential constraints and first order constraints are equivalent.
When $\bar{A}_{i}\bar{x}-\bar{b}_{i}\in\partial Q_{i}\backslash\{0\}$,
we may still be able to calculate the tangential constraints using
the material recalled in Section \ref{sec:var-analysis}.  

We illustrate the tangential problem with the example on second order
cone programming (SOCP).
\begin{example}
\label{exa:SOCP-easy}(Tangential problem in SOCP) Consider the SOCP
problem\begin{eqnarray*}
 &  & \min_{x}c^{T}x\\
 &  & \mbox{s.t. }\bar{A}_{i}x-\bar{b}_{i}\in Q_{k_{i}}\mbox{ for all }1\leq i\leq m,\end{eqnarray*}
where $Q_{d}\subset\mathbb{R}^{d}$ is the \emph{second order cone}
\[
Q_{d}:=\{w=(w_{0},\dots,w_{d-1})\in\mathbb{R}^{d}\mid\|(w_{1},\dots,w_{d-1})\|_{2}-w_{0}\leq0\}.\]
Given an optimal solution $\bar{x}$, we show how to obtain the tangential
constraint. If $\bar{A}_{i}\bar{x}-\bar{b}_{i}\in\intr(Q_{k_{i}})$,
then $T_{Q_{k_{i}}}(\bar{A}_{i}\bar{x}-\bar{b}_{i})=\mathbb{R}^{k_{i}}$
by Example \ref{exa:TFO}, and so the tangential constraint vanishes.
If $\bar{A}_{i}\bar{x}-\bar{b}_{i}=0$, then $T_{Q_{k_{i}}}(\bar{A}_{i}\bar{x}-\bar{b}_{i})=Q_{k_{i}}$
by Example \ref{exa:TFO}, so the tangential constraint is \[
\bar{A}_{i}\gamma+L_{i}(\Delta\mathcal{U}_{i})\subset Q_{k_{i}}.\]

\end{example}
We now consider the case $\bar{A}_{i}\bar{x}-\bar{b}_{i}\in\partial Q_{k_{i}}\backslash\{0\}$.
Let $\bar{z}=\bar{A}_{i}\bar{x}-\bar{b}_{i}$. In this case, $\bar{z}_{0}=\|(\bar{z}_{1},\dots,\bar{z}_{k_{i}-1})\|_{2}$.
The gradient of the map $(w_{0},w_{1},\dots,w_{k_{i}-1})\mapsto\|(w_{1},\dots,w_{k_{i}-1})\|_{2}-w_{0}$
at $\bar{z}$ is $(-1,\frac{(\bar{z}_{1},\dots,\bar{z}_{k_{i}-1})}{\|(\bar{z}_{1},\dots,\bar{z}_{k_{i}-1})\|_{2}})$.
Let $R:\mathbb{R}^{k_{i}}\to\mathbb{R}^{k_{i}}$ be the reflection
map \[
R(w_{0},w_{1},\dots,w_{k_{i}-1}):=(-w_{0},w_{1},\dots,w_{k_{i}-1}),\]
i.e., $R$ multiplies the $0$th coordinate by $-1$. The gradient
at $\bar{z}$ can also be written as $\frac{1}{\bar{z}_{0}}R\bar{z}$.
Therefore, by Proposition \ref{pro:tangent-examples},  \[
T_{Q_{k_{i}}}(\bar{A}_{i}\bar{x}-\bar{b}_{i})=\{w\in\mathbb{R}^{k_{i}}\mid w^{T}R\bar{z}\leq0\}.\]
Therefore, the tangential constraint is\[
\bar{A}_{i}\gamma+[(\Delta A_{i})\bar{x}-\Delta b_{i}]\in T_{Q_{k_{i}}}(\bar{A}_{i}\bar{x}-\bar{b}_{i})\mbox{ for all }(\Delta A_{i},\Delta b_{i})\in\Delta\mathcal{U}_{i}.\]
This can be written equivalently as \begin{eqnarray*}
 &  & \bar{z}^{T}R\bar{A}_{i}\gamma+\bar{z}^{T}R[(\Delta A_{i})\bar{x}-\Delta b_{i}]\leq0\mbox{ for all }(\Delta A_{i},\Delta b_{i})\in\Delta\mathcal{U}_{i},\\
 & \mbox{or } & \bar{z}^{T}R\bar{A}_{i}\gamma\leq-\max_{(\Delta A_{i},\Delta b_{i})\in\Delta\mathcal{U}_{i}}\bar{z}^{T}R[(\Delta A_{i})\bar{x}-\Delta b_{i}].\end{eqnarray*}

\section{\label{sec:Main-result}Main result: Approximation using the tangential
problem}

In Theorem \ref{thm:TFO-approx} we prove that if the uncertainty
set in a robust optimization problem is dilated or expanded, then
the robust optimal solution can be predicted from the exact solution
of the nonrobust problem and the tangential problem. 

We now prove a lemma needed for the proof of our main result. 
\begin{lem}
\label{lem:compex-sets-in-intr} (Compact sets in convex cones) Let
$D\subset\mathbb{R}^{n}$ be Clarke regular at $0$ and $C\subset\mathbb{R}^{n}$
be a compact convex set such that $\intr(T_{D}(0))\neq\emptyset$
and $C\subset T_{D}(0)$. Let $v\in\intr(T_{D}(0))$. Then for all
sufficiently small $\delta>0$, there exists $\bar{\epsilon}>0$ such
that $C+\delta v+\delta^{2}\mathbb{B}\subset\frac{1}{\epsilon}D$
for all $\epsilon\in(0,\bar{\epsilon}]$.\end{lem}
\begin{proof}
Since $v\in\intr(T_{D}(0))$, $v+\delta\mathbb{B}\subset\intr(T_{D}(0))$
for all sufficiently small $\delta>0$, and therefore $C+\delta v+\delta^{2}\mathbb{B}\subset\intr(T_{D}(0))$.
For every point $w\in C+\delta v+\delta^{2}\mathbb{B}$, we can find
a convex polyhedral set $P_{w}$ such that $\intr(P_{w})\neq\emptyset$
and $P_{w}\subset\intr(T_{D}(0))$. A compactness argument shows that
the set $C+\delta v+\delta^{2}\mathbb{B}$ is contained in the interior
of finitely many of these convex polyhedral sets, so there is a convex
polyhedral set $P$ such that $C+\delta v+\delta^{2}\mathbb{B}\subset P\subset\intr(T_{D}(0))$.

By the recession properties of tangent cones and the Clarke regularity
of $D$, there is an $\bar{\epsilon}>0$ such that $\epsilon\mbox{conv}(\{0\}\cup P)\subset D$
for all $\epsilon\in[0,\bar{\epsilon}]$ (see \cite[Exercise 6.34(a)]{RW98}.
The roots of this result on local recession vectors can be traced
back to \cite{R79}.). Therefore $C+\delta v+\delta^{2}\mathbb{B}\subset\mbox{conv}(\{0\}\cup P)\subset\frac{1}{\epsilon}D$
for all $\epsilon\in[0,\bar{\epsilon}]$ as needed.
\end{proof}
We also need material in set-valued analysis as presented in \cite[Chapters 4 and 5]{RW98}
for the proof of Theorem \ref{thm:TFO-approx}.
\begin{defn}
\cite[Definition 5.4]{RW98} (Set-valued continuity) We say that $S$
is a \emph{set-valued map}, denoted by $S:\mathbb{R}^{n}\rightrightarrows\mathbb{R}^{m}$,
if $S(x)\subset\mathbb{R}^{m}.$ A set-valued map $S$ is \emph{outer
semicontinuous (osc)} at $\bar{x}$ if \[
\limsup_{x\to\bar{x}}S(x)\subset S(\bar{x}),\]
or equivalently $\limsup_{x\to\bar{x}}S(x)=S(\bar{x})$, but \emph{inner
semicontinuous (isc)} at $\bar{x}$ if \[
\liminf_{x\to\bar{x}}S(x)\supset S(\bar{x}),\]
or equivalently when $S$ is closed-valued, $\liminf_{x\to\bar{x}}S(x)=S(\bar{x})$.
It is called \emph{continuous} at $\bar{x}$ if both conditions hold,
i.e., if $S(x)\to S(\bar{x})$ as $x\to\bar{x}$. Here, the \emph{outer
limit }$\limsup_{x\to\bar{x}}S(x)$ and the \emph{inner limit} $\liminf_{x\to\bar{x}}S(x)$
are defined by \begin{eqnarray*}
\limsup_{x\to\bar{x}}S(x) & := & \left\{ u\mid\exists x_{i}\to\bar{x},\exists u_{i}\to u\mbox{ with }u_{i}\in S(x_{i})\right\} \\
\liminf_{x\to\bar{x}}S(x) & := & \big\{u\mid\forall x_{i}\to\bar{x},\exists\{u_{i}\}\mbox{ with }u_{i}\in S(x_{i})\\
 &  & \quad\mbox{ s.t. }u\mbox{ is the limit of a subsequence of }\{u_{i}\}\big\}.\end{eqnarray*}
If $S$ maps to compact sets, continuity as defined by inner and outer
limits above is equivalent to continuity in the Pompieu-Hausdorff
distance, which is a metric in the subset of compact sets. We refer
to \cite{RW98} for more details. We also need to recall the definition
of epi-convergence. A sequence of functions $h_{i}:\mathbb{R}^{n}\to\mathbb{R}$
is said to \emph{epi-converge }to a function $h:\mathbb{R}^{n}\to\mathbb{R}$,
written $h_{i}\xrightarrow{e}h$, if $\epi(h_{i})\to\epi(h)$. The
history of epi-convergence can be traced back to the 1960's, and the
result we need for our proof (\cite[Theorem 7.33]{RW98}) can be traced
back to Salinetti (unpublished, but reported in \cite{RW84}) and
\cite{AW81}. See \cite[Chapter 7]{RW98}.
\end{defn}
Here is our theorem on the approximation properties of the tangential
problem. 
\begin{thm}
\label{thm:TFO-approx}(Approximation properties of Tangential problem)
Consider the robust optimization problem\begin{eqnarray}
 &  & \min_{x}\, f(x)\nonumber \\
 &  & \mbox{s.t. }[\bar{A}_{i}+\Delta A_{i}]x-[\bar{b}_{i}+\Delta b_{i}]\in Q_{i}\mbox{ for all }(\Delta A_{i},\Delta b_{i})\in\epsilon\Delta\mathcal{U}_{i}\mbox{ for all }i,\label{eq:RO_approx_1}\end{eqnarray}
and the tangential problem \begin{eqnarray}
 &  & \min_{\gamma}\, df(\bar{x})(\gamma)\nonumber \\
 &  & \mbox{s.t. }\bar{A}_{i}\gamma+[(\Delta A_{i})\bar{x}-\Delta b_{i}]\in T_{Q_{i}}(\bar{A}_{i}\bar{x}-\bar{b}_{i})\label{eq:FO_approx_1}\\
 &  & \mbox{ for all }(\Delta A_{i},\Delta b_{i})\in\Delta\mathcal{U}_{i}\mbox{ for all }i.\nonumber \end{eqnarray}
Let $\bar{\Gamma}$ be the set of optimal solutions to \eqref{eq:FO_approx_1},
$\Phi=\{x\mid\bar{A}x-\bar{b}\in Q\}$, and $\bar{x}$ be a solution
of the nominal problem $\min\{f(x)\mid\bar{A}x-\bar{b}\in Q\}$. Suppose 
\begin{enumerate}
\item $Q_{i}$ are closed sets that are Clarke regular at $\bar{A}_{i}\bar{x}-\bar{b}_{i}$, 
\item $\Delta\mathcal{U}_{i}$ are compact convex sets
\item $\bar{\Gamma}$ is bounded.
\item \label{enu:CQ}There is some $\gamma^{\prime}$ such that $\bar{A}_{i}\gamma^{\prime}\in\intr(T_{Q_{i}}(\bar{A}_{i}\bar{x}-\bar{b}_{i}))$
for all $i$. 
\item (Compactness) If $\{x_{i}\}\subset\Phi$ and $f(x_{i})\to f(\bar{x})$,
then $x_{i}\to\bar{x}$.
\item $f$ is locally Lipschitz and subdifferentially regular at $\bar{x}$.
\end{enumerate}
Then the set of cluster points of any sequence $\{\frac{1}{\epsilon}(\bar{x}_{\epsilon}-\bar{x})\}$,
where $\bar{x}_{\epsilon}$ is an optimal solution to \eqref{eq:RO_approx_1}
and $\epsilon\searrow0$, is a subset of $\bar{\Gamma}$. The objective
value of \eqref{eq:RO_approx_1} , say $\bar{v}_{\epsilon}$, has
an approximation $\bar{v}_{\epsilon}=\bar{v}+\epsilon\tilde{v}+o(\epsilon)$,
where $\tilde{v}$ is the objective value of \eqref{eq:FO_approx_1}.

In particular, if $\bar{\Gamma}$ contains only one element, say $\bar{\gamma}$,
then $\lim_{\epsilon\searrow0}\frac{1}{\epsilon}(\bar{x}_{\epsilon}-\bar{x})=\bar{\gamma}$,
or $\bar{x}_{\epsilon}\in\bar{x}+\epsilon\bar{\gamma}+o(\epsilon)$.\end{thm}
\begin{proof}
The proof of this result is broken up into four steps. In steps 1
to 3, we prove this result for the affine function $f(x)=c^{T}x+d$,
and $df(\bar{x})(\gamma)=c^{T}x$. In step 4, we use the observation
in Example \ref{exa:nonlin-obj} to treat the case where $f$ is locally
Lipschitz and subdifferentially regular at $\bar{x}$.

\textbf{Step 1: Rewriting the robust optimization problem \eqref{eq:RO_approx_1}.}

We rewrite the constraint in the robust optimization problem.\begin{eqnarray*}
 &  & [\bar{A}_{i}+\Delta A_{i}]x-[\bar{b}_{i}+\Delta b_{i}]\in Q_{i}\\
 & \Leftrightarrow & [\bar{A}_{i}+\Delta A_{i}]x-\bar{A}_{i}\bar{x}-\Delta b_{i}\in Q_{i}-[\bar{A}_{i}\bar{x}-\bar{b}_{i}]\\
 & \Leftrightarrow & \bar{A}_{i}(x-\bar{x})+[(\Delta A_{i})\bar{x}-\Delta b_{i}]+(\Delta A_{i})(x-\bar{x})\in Q_{i}-[\bar{A}_{i}\bar{x}-\bar{b}_{i}].\end{eqnarray*}
Hence, \begin{eqnarray*}
 &  & A_{i}x-b_{i}\in Q_{i}\mbox{ for all }(\Delta A_{i},\Delta b_{i})\in\epsilon\Delta\mathcal{U}_{i}\\
 & \Leftrightarrow & \bar{A}_{i}(x-\bar{x})+[(\Delta A_{i})\bar{x}-\Delta b_{i}]+(\Delta A_{i})(x-\bar{x})\in Q_{i}-[\bar{A}_{i}\bar{x}-\bar{b}_{i}]\\
 &  & \quad\quad\mbox{ for all }(\Delta A_{i},\Delta b_{i})\in\epsilon\Delta\mathcal{U}_{i}.\end{eqnarray*}
The next step is to scale the variables $\Delta A_{i}$ and $\Delta b_{i}$
so that the $\epsilon$ vanishes from the expression $\epsilon\Delta\mathcal{U}_{i}$.
This gives\begin{eqnarray}
 &  & \bar{A}_{i}(x-\bar{x})+[(\Delta A_{i})\bar{x}-\Delta b_{i}]+(\Delta A_{i})(x-\bar{x})\in Q_{i}-[\bar{A}_{i}\bar{x}-\bar{b}_{i}]\nonumber \\
 &  & \quad\quad\mbox{ for all }(\Delta A_{i},\Delta b_{i})\in\epsilon\Delta\mathcal{U}_{i}\nonumber \\
 & \Leftrightarrow & \bar{A}_{i}\Big[\frac{1}{\epsilon}(x-\bar{x})\Big]+[(\Delta A_{i})\bar{x}-\Delta b_{i}]+\epsilon(\Delta A_{i})\Big[\frac{1}{\epsilon}(x-\bar{x})\Big]\nonumber \\
 &  & \quad\quad\in\frac{1}{\epsilon}\big[Q_{i}-[\bar{A}_{i}\bar{x}-\bar{b}_{i}]\big]\mbox{ for all }(\Delta A_{i},\Delta b_{i})\in\Delta\mathcal{U}_{i}\nonumber \\
 & \Leftrightarrow & \bar{A}_{i}\gamma_{\epsilon}+[(\Delta A_{i})\bar{x}-\Delta b_{i}]+\epsilon(\Delta A_{i})\gamma_{\epsilon}\label{eq:rob_approx_1}\\
 &  & \quad\quad\in\frac{1}{\epsilon}\big[Q_{i}-[\bar{A}_{i}\bar{x}-\bar{b}_{i}]\big]\mbox{ for all }(\Delta A_{i},\Delta b_{i})\in\Delta\mathcal{U}_{i},\nonumber \end{eqnarray}
where $\gamma_{\epsilon}:=\frac{1}{\epsilon}(x-\bar{x})$ in the final
expression. We see that as $\epsilon\searrow0$, the expressions in
\eqref{eq:rob_approx_1} converge to the corresponding expressions
for the tangential constraints. 

Let $\Gamma_{\epsilon}$ denote the set of all feasible $\gamma_{\epsilon}$
for the robust problem with parameter $\epsilon$, and $\Gamma$ denote
the set of all feasible $\gamma$ for the tangential problem. Similarly,
let $\bar{\Gamma}_{\epsilon}$ and $\bar{\Gamma}$ denote the set
of optimal solutions to the corresponding problems. It is elementary
to check that the sets $\Gamma_{\epsilon}$, $\Gamma$, $\bar{\Gamma}_{\epsilon}$
and $\bar{\Gamma}$ are all closed.

\textbf{Step 2: $\lim_{\epsilon\searrow0}\Gamma_{\epsilon}=\Gamma$.}

Suppose that $\{\gamma_{j}\}$ is a sequence of feasible solutions
to the robust problem with parameter $\epsilon_{j}$, that is $\gamma_{j}\in\Gamma_{\epsilon_{j}}$.
Then each $\gamma_{j}$ satisfies the formula in \eqref{eq:rob_approx_1}
with parameter $\epsilon_{j}$. It is clear that any limit of $\{\gamma_{j}\}$
is a feasible solution of the tangential problem \eqref{eq:FO_approx_1},
so $\limsup_{\epsilon\searrow0}\Gamma_{\epsilon}\subset\Gamma$.

Next, we show that $\liminf_{\epsilon\searrow0}\Gamma_{\epsilon}\supset\Gamma$.
Suppose $\tilde{\gamma}\in\Gamma$. We need to show that for any
choice of $\epsilon_{j}\searrow0$, we can find $\gamma_{j}\in\Gamma_{\epsilon_{j}}$
such that $\tilde{\gamma}=\lim_{j\to\infty}\gamma_{j}$. Recall that
$\tilde{\gamma}$ satisfies \[
\bar{A}_{i}\tilde{\gamma}+L_{i}(\Delta\mathcal{U}_{i})\subset T_{Q_{i}}(\bar{A}_{i}\bar{x}-\bar{b}_{i}),\]
where the linear map $L_{i}:\mathbb{R}^{m_{i}\times n}\times\mathbb{R}^{m_{i}}\to\mathbb{R}^{m_{i}}$
by $L_{i}(\Delta A_{i},\Delta b_{i})=(\Delta A_{i})\bar{x}-\Delta b_{i}$.
It follows from the convexity of $\Delta\mathcal{U}_{i}$ that $L_{i}(\Delta\mathcal{U}_{i})$
is convex. Since $\bar{A}_{i}\gamma^{\prime}\in\intr(T_{Q_{i}}(\bar{A}_{i}\bar{x}-\bar{b}_{i}))$,
we can apply Lemma \ref{lem:compex-sets-in-intr} to tell us that
for all sufficiently small $\delta>0$, there is some $\bar{\epsilon}>0$
such that $[\bar{A}_{i}(\tilde{\gamma}+\delta\gamma^{\prime})+L_{i}(\Delta\mathcal{U}_{i})+\delta^{2}\mathbb{B}]\subset\frac{1}{\epsilon}[Q_{i}-(\bar{A}_{i}\bar{x}-\bar{b}_{i})]$
for all $\epsilon\in[0,\bar{\epsilon}]$. Therefore \[
\bar{A}_{i}(\tilde{\gamma}+\delta\gamma^{\prime})+L_{i}(\Delta\mathcal{U}_{i})+\delta^{2}\mathbb{B}\subset\frac{1}{\epsilon}\big[Q_{i}-[\bar{A}_{i}\bar{x}-\bar{b}_{i}]\big].\]
If $\epsilon<\bar{\epsilon}$ and $\epsilon\max\{\|\Delta A_{i}\|\mid(\Delta A_{i},\Delta b_{i})\in\Delta\mathcal{U}_{i}\}\|(\tilde{\gamma}+\delta\gamma^{\prime})\|<\delta^{2}$
, then for all $(\Delta A_{i},\Delta b_{i})\in\Delta\mathcal{U}_{i}$,
\begin{eqnarray}
 &  & \bar{A}_{i}(\tilde{\gamma}+\delta\gamma^{\prime})+[(\Delta A_{i})\bar{x}-\Delta b_{i}]+\epsilon(\Delta A_{i})(\tilde{\gamma}+\delta\gamma^{\prime})\label{eq:Aw-chain}\\
 & \in & \bar{A}_{i}(\tilde{\gamma}+\delta\gamma^{\prime})+L_{i}(\Delta\mathcal{U}_{i})+\epsilon\|\Delta A_{i}\|\|(\tilde{\gamma}+\delta\gamma^{\prime})\|\mathbb{B}\nonumber \\
 & \subset & \bar{A}_{i}(\tilde{\gamma}+\delta\gamma^{\prime})+L_{i}(\Delta\mathcal{U}_{i})+\delta^{2}\mathbb{B}\nonumber \\
 & \subset & \frac{1}{\epsilon}\big[Q_{i}-[\bar{A}_{i}\bar{x}-\bar{b}_{i}]\big].\nonumber \end{eqnarray}
With this observation, we can choose a sequence $\delta_{j}\searrow0$
such that $(\tilde{\gamma}+\delta_{j}\gamma^{\prime})\in\Gamma_{\epsilon_{j}}$,
which gives $\liminf_{\epsilon\searrow0}\Gamma_{\epsilon}\supset\Gamma$
as needed. 

\textbf{Step 3: $\limsup_{\epsilon\searrow0}\bar{\Gamma}_{\epsilon}\subset\bar{\Gamma}$.}

Recall that for a closed set $D\subset\mathbb{R}^{n}$, the \emph{indicator
function} $\delta_{D}:\mathbb{R}^{n}\to\mathbb{R}\cup\{\infty\}$
is defined by \[
\delta_{D}(x):=\begin{cases}
0 & \mbox{if }x\in D\\
\infty & \mbox{otherwise.}\end{cases}\]
Define $h_{\epsilon}:\mathbb{R}^{n}\to\mathbb{R}$ and $h:\mathbb{R}^{n}\to\mathbb{R}$
by $h_{\epsilon}(\gamma):=c^{T}\gamma+\delta_{\Gamma_{\epsilon}}(\gamma)$
and $h(\gamma):=c^{T}\gamma+\delta_{\Gamma}(\gamma)$. Since $\Gamma_{\epsilon}\to\Gamma$,
we have $\epi(\delta_{\Gamma_{\epsilon}})=\Gamma_{\epsilon}\times[0,\infty)\to\Gamma\times[0,\infty)=\epi(\delta_{\Gamma})$
(by \cite[Exercise 4.29(a)]{RW98}), so in other words $\delta_{\Gamma_{\epsilon}}\xrightarrow{e}\delta_{\Gamma}$.
By \cite[Exercise 7.8(a)]{RW98} we have $h_{\epsilon}\xrightarrow{e}h$.
We seek to apply \cite[Theorem 7.33]{RW98}, which gives us the result
we need. Before we can do so, we have to check that $\{h_{\epsilon}\}$
is eventually level bounded, that is, for any $\alpha\in\mathbb{R}$,
we have $\cup_{[0,\epsilon]}\{\gamma\mid h_{\epsilon}(\gamma)\leq\alpha\}$
being bounded for some $\epsilon>0$. 

Recall $\Phi=\{x\mid\bar{A}x-\bar{b}\in Q\}$. By Proposition \ref{pro:Condn-Phi},
the boundedness of $\bar{\Gamma}$ is equivalent to $\{c\}^{\perp}\cap T_{\Phi}(\bar{x})=\{0\}$.
Suppose on the contrary that $\{h_{\epsilon}\}$ is not eventually
level bounded. Then there is some $\alpha$ and sequences $\{\epsilon_{i}\}$
and $\{\gamma_{i}\}$ such that $\epsilon_{i}\searrow0$, $\{\gamma_{i}\}$
is unbounded, and $h_{\epsilon_{i}}(\gamma_{i})\leq\alpha$ (or equivalently,
$\gamma_{i}\in\Gamma_{\epsilon_{i}}$ and $c^{T}\gamma_{i}\leq\alpha$).
Let us write $x_{i}=\bar{x}+\epsilon_{i}\gamma_{i}$. Since $x_{i}\in\Phi$
and \[
c^{T}\bar{x}\leq c^{T}x_{i}=c^{T}\bar{x}+\epsilon_{i}c^{T}\gamma_{i}\leq c^{T}\bar{x}+\epsilon_{i}\alpha,\]
we have $c^{T}x_{i}\to c^{T}\bar{x}$. Note that this also gives us
$\alpha\geq0$. By our compactness assumption, $x_{i}\to\bar{x}$,
which means that $|\epsilon_{i}\gamma_{i}|\to0$. This means that
$\frac{\gamma_{i}}{|\gamma_{i}|}$ converges to a vector $\gamma_{\infty}$
in $T_{\Phi}(0)\backslash\{0\}$. Observe that $c^{T}\gamma_{\infty}=\lim_{i\to\infty}c^{T}\frac{\gamma_{i}}{|\gamma_{i}|}\leq\lim_{i\to\infty}\frac{\alpha}{|\gamma_{i}|}=0$,
and that $c^{T}\gamma_{i}\geq0$, so $c^{T}\gamma_{\infty}=0$. This
is a contradiction to $\{c\}^{\perp}\cap T_{\Phi}(\bar{x})=\{0\}$.
We can thus apply \cite[Theorem 7.33]{RW98} to conclude that $\min\{c^{T}\gamma\mid\gamma\in\Gamma_{\epsilon}\}$
converges to $\min\{c^{T}\gamma\mid\gamma\in\Gamma\}$ and $\limsup_{\epsilon\searrow0}\bar{\Gamma}_{\epsilon}\subset\bar{\Gamma}$,
ending the proof of the theorem for the linear case.

\textbf{Step 4: Locally Lipschitz subdifferentially regular $f$ at
$\bar{x}$.}

Consider the problem \begin{eqnarray}
 &  & \min_{x,t}\, t\nonumber \\
 &  & \mbox{s.t. }[\bar{A}_{i}+\Delta A_{i}]x-[\bar{b}_{i}+\Delta b_{i}]\in Q_{i}\mbox{ for all }(\Delta A_{i},\Delta b_{i})\in\epsilon\Delta\mathcal{U}_{i}\mbox{ for all }i,\label{eq:nonlin-robust}\\
 &  & \mbox{and }(x,t)\in\epi(f)=\{(x,t)\mid f(x)\leq t\}.\nonumber \end{eqnarray}
and the tangential problem \begin{eqnarray}
 &  & \min_{\gamma,s}\, s\nonumber \\
 &  & \mbox{s.t. }\bar{A}_{i}\gamma+[(\Delta A_{i})\bar{x}-\Delta b_{i}]\in T_{Q_{i}}(\bar{A}_{i}\bar{x}-\bar{b}_{i})\label{eq:nonlin-tangent}\\
 &  & \phantom{\mbox{s.t. }}\mbox{ for all }(\Delta A_{i},\Delta b_{i})\in\Delta\mathcal{U}_{i}\mbox{ for all }i.\nonumber \\
 &  & \mbox{and }(\gamma,s)\in\Tepi(\bar{x},f(\bar{x})).\nonumber \end{eqnarray}

The robust and tangential problems are equivalent to the respective
problems \eqref{eq:RO_approx_1} and \eqref{eq:FO_approx_1} in the
statement of the theorem. We now show that if conditions (1) to (5)
in the theorem statement are satisfied for \eqref{eq:RO_approx_1}
and \eqref{eq:FO_approx_1}, then these conditions hold for \eqref{eq:nonlin-robust}
and \eqref{eq:nonlin-tangent} as well.

For condition (1), we need only to check that $\epi(f)$ is Clarke
regular at $(\bar{x},f(\bar{x}))$, which is immediate from the subdifferential
regularity of $f$ at $\bar{x}$. Condition (2) is straightforward.
For condition (3), we note that the set of minimizers of \eqref{eq:nonlin-tangent}
is just $\bar{\Gamma}\times\{df(\bar{x})(\hat{\gamma})\}$, where
$\hat{\gamma}$ is any element in $\bar{\gamma}$. The set $\bar{\Gamma}\times\{df(\bar{x})(\hat{\gamma})\}$
is bounded if and only if $\bar{\Gamma}$ is bounded. 

We further assume that $f$ is locally Lipschitz at $\bar{x}$ with
Lipschitz modulus $\kappa$. For condition (4), suppose $\gamma^{\prime}$
is a vector such that $\bar{A}_{i}\gamma^{\prime}\in\intr(T_{Q_{i}}(\bar{A}_{i}\bar{x}-\bar{b}_{i}))$
for all $i$. Since $(\mathbf{0},1)\in\intr(\Tepi(\bar{x},f(\bar{x})))$,
we have $(\gamma^{\prime},(\kappa+1)|\gamma^{\prime}|)\in\intr(\Tepi(\bar{x},f(\bar{x})))$,
which verifies condition (4). 

We also need to check that given $f(x_{i})\to f(\bar{x})$ and $\{x_{i}\}\subset\Phi$
implies $x_{i}\to\bar{x}$, we have the compactness condition that
$t_{i}\to f(\bar{x})$, $(x_{i},t_{i})$ satisfies $f(x_{i})\leq t_{i}$
and $\{x_{i}\}\subset\Phi$ implies $(x_{i},t_{i})\to(\bar{x},f(\bar{x}))$.
Since $f(\bar{x})\leq f(x_{i})\leq t_{i}$ and $t_{i}\to f(\bar{x})$,
we have $f(x_{i})\to f(\bar{x})$, which gives $x_{i}\to\bar{x}$,
and thus $(x_{i},t_{i})\to(\bar{x},f(\bar{x}))$ as needed.
\end{proof}
We now take a closer look at step 4 of the proof of Theorem \ref{thm:TFO-approx}.
Consider the general case where $\bar{\Gamma}_{\epsilon}=\arg\min\{c^{T}\gamma\mid\gamma\in\Gamma_{\epsilon}\}$,
$\bar{\Gamma}=\arg\min\{c^{T}\gamma\mid\gamma\in\Gamma\}$ and $\Gamma=\lim_{\epsilon\searrow0}\Gamma_{\epsilon}$.
It may turn out that $\limsup_{\epsilon\searrow0}\bar{\Gamma}_{\epsilon}\subsetneq\bar{\Gamma}$,
as the example in Figure \ref{fig:V_e_V} shows. Example \ref{exa:square-eg}
shows that it is possible for $\bar{\Gamma}$ to be bounded but not
be a singleton set. In such cases, it is possible that $\bar{\Gamma}_{\epsilon}$
is a singleton set, which occurs when the function $f$ is strictly
convex for example. 

\begin{figure}[h]
\includegraphics[scale=0.5]{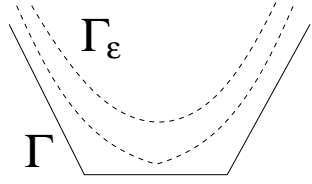}

\caption{\label{fig:V_e_V}$\lim_{\epsilon\to0}\Gamma_{\epsilon}=\Gamma$,
but $\limsup_{\epsilon\to0}\bar{\Gamma}_{\epsilon}\subsetneq\bar{\Gamma}$.}

\end{figure}

We make an observation on the condition $\bar{A}_{i}\gamma^{\prime}\in\intr(T_{Q_{i}}(\bar{A}_{i}\bar{x}-\bar{b}_{i}))$
for all $i$ in Theorem \ref{thm:TFO-approx}.
\begin{example}
(Constraint qualification in tangential problem) The optimization
problem\begin{eqnarray*}
 &  & \min_{x\in\mathbb{R}^{2}}x_{2}\\
 &  & \mbox{s.t. }x\in Q_{1}:=\{x\mid-x_{1}+x_{2}^{2}\leq0\}\\
 &  & \phantom{\mbox{s.t. }}x\in Q_{2}:=\{x\mid x_{1}+x_{2}^{2}\leq0\}\end{eqnarray*}
can be written equivalently as \begin{eqnarray*}
 &  & \min_{x\in\mathbb{R}^{2}}x_{2}\\
 &  & \mbox{s.t. }x\in Q_{1}\cap Q_{2}=\{(0,0)\},\end{eqnarray*}
and the solution for both problems is $\bar{x}=(0,0)$. The tangential
approximation for the first problem is \begin{eqnarray*}
 &  & \min_{\gamma\in\mathbb{R}^{2}}\gamma_{2}\\
 &  & \left.\begin{array}{c}
\mbox{s.t. }\gamma\in T_{Q_{1}}(\bar{x})=\mathbb{R}_{+}\times\mathbb{R}\\
\phantom{\mbox{s.t. }}\gamma\in T_{Q_{2}}(\bar{x})=\mathbb{R}_{-}\times\mathbb{R}\end{array}\right\} \implies\gamma\in\{0\}\times\mathbb{R},\end{eqnarray*}
while the tangential approximation for the second problem is \begin{eqnarray*}
 &  & \min_{\gamma\in\mathbb{R}^{2}}\gamma_{2}\\
 &  & \mbox{s.t. }\gamma\in T_{Q_{1}\cap Q_{2}}(\bar{x})=\{(0,0)\}.\end{eqnarray*}
Clearly, the solutions for the two problems are different. This example
shows that depending on how the optimization problem is written, the
tangential problems may not be equivalent and may have different solutions.
But note that in this case, $T_{Q_{1}\cap Q_{2}}(\bar{x})\subsetneq T_{Q_{1}}(\bar{x})\cap T_{Q_{2}}(\bar{x})$,
which implies that the MFCQ does not hold, which in turn implies that
there is no vector $\gamma^{\prime}$ such that $\bar{A}_{i}\gamma^{\prime}\in\intr(T_{Q_{i}}(\bar{A}_{i}\bar{x}-\bar{b}_{i}))$
for all $i$. The condition $\bar{A}_{i}\gamma^{\prime}\in\intr(T_{Q_{i}}(\bar{A}_{i}\bar{x}-\bar{b}_{i}))$
for all $i$ in Theorem \ref{thm:TFO-approx} ensures that the MFCQ
holds whenever a group of the $\{(A_{i},b_{i})\}_{i=1}^{m}$ contains
repetitions, ensuring that the tangent cones of the intersections
is the intersections of the corresponding tangent cones through Proposition
\ref{pro:tangent-NMFCQ}. 

Theorem \ref{thm:TFO-approx} shows that the decrease in objective
function of the robust optimization problem is differentiable in the
size of the uncertainty set at $\epsilon=0$. This observation can
help give an approximate of the maximum robustness one can afford
if the objective is to be above a certain value. The tangential problem
also shows that the variables that we should estimate or measure more
accurately are those which make the set $L(\Delta\mathcal{U})=\{(\Delta A)\bar{x}-\Delta b\mid(\Delta A,\Delta b)\in\Delta\mathcal{U}\}$
small. For example, if $\bar{x}=0$, then more effort should be spent
on determining $\bar{b}$ accurately rather than entries in $\bar{A}$.
Likewise, $\bar{x}$ determines which variables in $\bar{A}$ should
be measured more accurately than others.

Besides these analytical properties, Theorem \ref{thm:TFO-approx}
shows that solving the tangential problem can help to obtain a good
approximate of the robust solution. The robust optimization problem
is known to be more computationally expensive than the original problem,
so it will take more effort to obtain a desired level of accuracy.
With the tangential problem, we can make use of the previously calculated
optimization problem to obtain an approximate of the robust solution.
Such an approximation is likely to be simpler than the robust optimization
problem (for example, for LP in \eqref{eq:LP-easy} and for SOCP in
Example \ref{exa:SOCP-easy}, though it still may be computationally
difficult), and need not be computed to very high accuracy to obtain
a good approximate of the robust optimization problem.\end{example}
\begin{rem}
(Relaxing the constraint qualification in Theorem \ref{thm:TFO-approx})
The existence of $\gamma^{\prime}$ such that $\bar{A}\gamma^{\prime}\in\intr(T_{Q}(\bar{A}\bar{x}-\bar{b}))$
can be relaxed slightly if more structure is known about the set $Q$.
All we need is for the chain of inclusions \eqref{eq:Aw-chain} to
hold. For example, if $Q_{i}$ is polyhedral and there is some vector
$\gamma^{\prime}$ such that \[
\bar{A}_{i}\gamma^{\prime}\in\rint\big(\spanv(\Pi_{i}(\Delta\mathcal{U}_{i}))\big)\cap T_{Q_{i}}(\bar{A}_{i}\bar{x}-\bar{b}_{i}),\]
where $\Pi_{i}(\Delta\mathcal{U}_{i})=\{\Delta A_{i}\mid(\Delta A_{i},\Delta b_{i})\in\Delta\mathcal{U}_{i}\mbox{ for some }\Delta b_{i}\}$,
then the chain of inclusions \eqref{eq:Aw-chain} would hold as well.
A finer analysis on $\Pi_{i}(\Delta\mathcal{U})$ and local recession
vectors can give stronger results.
\end{rem}

\section{\label{sec:First-order}First order optimality conditions of the
tangential problem}

In this section, we discuss first order optimality conditions of the
tangential problem, which can be useful for designing specialized
numerical methods for the tangential  problem. In view of Theorem
\ref{thm:TFO-approx}, we also give sufficient conditions for $\bar{\Gamma}$
to be bounded and for $\bar{\Gamma}$ to be a singleton. 

As explained in \cite[Chapters 5-8]{BGN09}, the robust optimization
problem is computationally tractable if either $Q$ is polyhedral
or $\Delta\mathcal{U}$ is polyhedral, while most other problems encountered
in practice are not computationally tractable. Recall that a typical
constraint in a robust optimization problem whose nominal solution
is $\bar{x}$ is \[
[\bar{A}+\Delta A](\bar{x}+\Delta x)-[\bar{b}+\Delta b]\in Q\mbox{ for all }(\Delta A,\Delta b)\in\Delta\mathcal{U}.\]
If there were no uncertainty in the matrix $\bar{A}$, then the constraint
can be written as\begin{equation}
\bar{A}(\Delta x)-\Delta b\in Q-[\bar{A}\bar{x}-\bar{b}]\mbox{ for all }(\Delta A,\Delta b)\in\Delta\mathcal{U}.\label{eq:robust-fixed-A}\end{equation}
Recall that the tangential constraint is of the form \[
\bar{A}\gamma+[(\Delta A)\bar{x}-\Delta b]\in T_{Q}(\bar{A}\bar{x}-\bar{b})\mbox{ for all }(\Delta A,\Delta b)\in\Delta\mathcal{U}.\]
If the variable $\Delta x$ in \eqref{eq:robust-fixed-A} were replaced
by $\gamma$, then we see that \eqref{eq:robust-fixed-A} is similar
to the tangential constraint. In other words, the uncertainty in $\bar{A}$
is transferred to $\bar{b}$ in the tangential constraint through
$(\Delta A)\bar{x}$. Tangential problems are still often hard to
compute efficiently, but the additional structure may be exploited
for designing specialized methods. The case where $T_{Q}(\bar{A}\bar{x}-\bar{b})$
is polyhedral is just robust linear programming and is easy, while
the tangential problem for $L(\Delta\mathcal{U})$ being polyhedral
reduces to optimizing over the cone $T_{Q}(\bar{A}\bar{x}-\bar{b})$,
as illustrated below.
\begin{example}
(Polyhedral $L(\Delta\mathcal{U})$) Suppose $L(\Delta\mathcal{U})$
is a polyhedral compact set of the form $L(\Delta\mathcal{U})=\conv(v_{1},\dots,v_{J})$.
The tangential problem \begin{eqnarray*}
 &  & \min_{\gamma}h(\gamma)\\
 &  & \mbox{s.t. }\bar{A}\gamma+L(\Delta\mathcal{U})\subset T_{Q}(\bar{A}\bar{x}-\bar{b}),\end{eqnarray*}
where $h$ is sublinear, is equivalent to \begin{eqnarray*}
 &  & \min_{\gamma}h(\gamma)\\
 &  & \mbox{s.t. }\bar{A}\gamma+v_{j}\in T_{Q}(\bar{A}\bar{x}-\bar{b})\mbox{ for all }1\leq j\leq J.\end{eqnarray*}

\end{example}
For this section, we define the sets $\Phi$ and $\Psi$ by \begin{eqnarray}
\Phi & := & \{x:\bar{A}x-\bar{b}\in Q\}.\label{eq:Phi}\\
\Psi & := & \{y:y+L(\Delta\mathcal{U})\subset T_{Q}(\bar{A}\bar{x}-\bar{b})\}.\nonumber \end{eqnarray}
Hence $\Gamma$ can be written similarly as \begin{eqnarray}
\Gamma & = & \{\gamma:\bar{A}\gamma+L(\Delta\mathcal{U})\subset T_{Q}(\bar{A}\bar{x}-\bar{b})\}\label{eq:Gamma}\\
 & = & \{\gamma:\bar{A}\gamma\in\Psi\}.\nonumber \end{eqnarray}
Recall also that \begin{equation}
\bar{\Gamma}=\arg\min\{df(\bar{x})(\gamma)\mid\gamma\in\Gamma\}.\label{eq:barGamma}\end{equation}

\begin{prop}
\label{pro:T_Phi}(Tangent space of feasible set) Suppose $Q$ is
Clarke regular at $\bar{A}\bar{x}-\bar{b}$. If there is a vector
$\gamma^{\prime}$ such that $\bar{A}\gamma^{\prime}\in\rint(T_{Q}(\bar{A}\bar{x}-\bar{b}))$,
then $T_{\Phi}(\bar{x})=\{\gamma:\bar{A}\gamma\in T_{Q}(\bar{A}\bar{x}-\bar{b})\}$,
where $\Phi$ is defined in \eqref{eq:Phi}.\end{prop}
\begin{proof}
This follows directly from \cite[Theorem 6.31]{RW98}. The constraint
qualification in that Theorem is equivalent to the condition in our
statement through \cite[Exercise 6.39(b)]{RW98}.
\end{proof}
The normal cone $N_{\Gamma}(\gamma)$ can be estimated from the image
$N_{\Psi}(\bar{A}\gamma)$.
\begin{prop}
\label{pro:normal-of-tangent-feasible}(Normal cone of tangent feasible
set) Suppose there is a vector $\gamma^{\prime}$ such that $\bar{A}\gamma^{\prime}\in\intr(T_{Q}(\bar{A}\bar{x}-\bar{b}))$,
and $Q$ is Clarke regular at $\bar{A}\bar{x}-\bar{b}$. Then the
normal cone $N_{\Gamma}(\gamma)$ equals $\{\bar{A}^{T}v\mid v\in N_{\Psi}(\bar{A}\gamma)\}$,
where $\Gamma$ and $\Psi$ are defined in \eqref{eq:Phi} and \eqref{eq:Gamma}.\end{prop}
\begin{proof}
Note that $\Gamma$ can be written in terms of $\Psi$ as $\Gamma=\{\gamma:\bar{A}\gamma\in\Psi\}$.
The result follows directly from \cite[Theorem 6.14]{RW98}, though
we still have to check the constraint qualification condition there.
Through \cite[Exercise 6.39(b)]{RW98}, the constraint qualification
condition required is that there is a vector $\gamma^{\prime}$ such
that $\bar{A}\gamma^{\prime}\in\intr(T_{\Psi}(\bar{A}\gamma))$. Note
that the recession cone of $\Psi$ is $T_{Q}(\bar{A}\bar{x}-\bar{b})$.
This means that $T_{Q}(\bar{A}\bar{x}-\bar{b})\subset T_{\Psi}(\bar{A}\gamma)$,
which shows that $\bar{A}\gamma^{\prime}\in\intr(T_{Q}(\bar{A}\bar{x}-\bar{b}))$
implies the constraint qualification condition. The conclusion is
straightforward.
\end{proof}
One notices that if there is a vector $\gamma^{\prime}$ such that
$\bar{A}\gamma^{\prime}\in\intr(T_{Q}(\bar{A}\bar{x}-\bar{b}))$ and
$0\in\intr(L(\Delta\mathcal{U}))$ (which holds when $0\in\intr(\Delta\mathcal{U})$),
then the tangential problem is feasible. 

If a vector $\bar{x}$ is a solution of the problem $\min\{f(x)\mid x\in D\}$,
where $f$ is locally Lipschitz and subdifferentially regular at $\bar{x}$
then it is well known that there is a $c\in\partial f(\bar{x})$ such
that $-c\in N_{D}(\bar{x})$ (see \cite[Theorem 8.15]{RW98} for example).
We prove the following lemmas, whose proofs do not seem easy to find.
\begin{lem}
\label{lem:strict-loc-min}(Strict minimizers) Let $D\subset\mathbb{R}^{n}$
be Clarke regular at $\bar{x}$, $f:\mathbb{R}^{n}\to\mathbb{R}$
be locally Lipschitz and subdifferentially regular at $\bar{x}$.
Then $-c\in\intr(N_{D}(\bar{x}))$ for some $c\in\partial f(\bar{x})$
implies that $\bar{x}$ is a strict local minimizer of $\min\{f(x)\mid x\in D\}$.\end{lem}
\begin{proof}
Since $-c\in\intr(N_{D}(\bar{x}))$, there is some $\delta>0$ such
that $-c+\delta\mathbb{B}\subset N_{D}(\bar{x})$, and so there is
a neighborhood $O$ of $\bar{x}$ such that for $x\in D\cap O$, \begin{eqnarray*}
\Big[-c+\delta\frac{x-\bar{x}}{|x-\bar{x}|}\Big]^{T}(x-\bar{x}) & \leq & \frac{\delta}{2}|x-\bar{x}|\\
\Rightarrow c^{T}(x-\bar{x}) & \geq & \frac{\delta}{2}|x-\bar{x}|.\end{eqnarray*}
By restricting $O$ if necessary, for $x\in D\cap O$, we have $f(x)\geq f(\bar{x})+df(\bar{x})(x-\bar{x})-\frac{\delta}{4}|x-\bar{x}|$.
The well known characterization of the subderivative in terms of support
functions (see for example \cite{Cla83,R70,RW98}) gives \begin{eqnarray*}
f(x) & \geq & f(\bar{x})+df(\bar{x})(x-\bar{x})-\frac{\delta}{4}|x-\bar{x}|\\
 & = & f(\bar{x})+\max_{v\in\partial f(\bar{x})}v^{T}(x-\bar{x})-\frac{\delta}{4}|x-\bar{x}|\\
 & \geq & f(\bar{x})+c^{T}(x-\bar{x})-\frac{\delta}{4}|x-\bar{x}|\\
 & \geq & f(\bar{x})+\frac{\delta}{4}|x-\bar{x}|.\end{eqnarray*}
Therefore $\bar{x}$ is a strict local minimizer.\end{proof}
\begin{lem}
\label{lem:equiv-strict}(Equivalence of strict minimizer condition)
Let $D$ be Clarke regular at $\bar{x}$. When $f$ is $\mathcal{C}^{1}$
at $\bar{x}$, $\nabla f(\bar{x})\neq0$ and $-\nabla f(\bar{x})\in N_{D}(\bar{x})$
(which holds when $\bar{x}$ is a local minimizer of $\min\{f(x)\mid x\in D\}$),
then the conditions $-\nabla f(\bar{x})\in\intr(N_{D}(\bar{x}))$
and $\{\nabla f(\bar{x})\}^{\perp}\cap T_{D}(\bar{x})=\{0\}$ are
equivalent.\end{lem}
\begin{proof}
The fact that $-\nabla f(\bar{x})\in N_{D}(\bar{x})$ when $\bar{x}$
is a local minimizer of $\min\{f(x)\mid x\in D\}$ is well known (see
\cite[Theorem 6.12]{RW98} for example). 

Suppose that $-\nabla f(\bar{x})\in\intr(N_{D}(\bar{x}))$, and $v\in\{\nabla f(\bar{x})\}^{\perp}\cap T_{D}(\bar{x})$.
Then $v\in T_{D}(\bar{x})=[N_{D}(\bar{x})]^{*}$ (the polar cone of
$N_{D}(\bar{x})$). In other words, $v^{T}s\leq0$ for all $s\in N_{D}(\bar{x})$.
If $v\neq0$, since $-\nabla f(\bar{x})\in\intr(N_{D}(\bar{x}))$,
there is some $\epsilon>0$ such that $\nabla f(\bar{x})+\epsilon v\in N_{D}(\bar{x})$,
which gives $v^{T}(\nabla f(\bar{x})+\epsilon v)\leq0$, and thus
$v^{T}\nabla f(\bar{x})<0$. Since $v\in\{\nabla f(\bar{x})\}^{\perp}$,
this means $v=0$, so $\{\nabla f(\bar{x})\}^{\perp}\cap T_{D}(\bar{x})=\{0\}$.

Next, suppose $-\nabla f(\bar{x})\in N_{D}(\bar{x})\backslash\intr(N_{D}(\bar{x}))$.
Then $N_{N_{D}(\bar{x})}(-\nabla f(\bar{x}))\supsetneq\{0\}$. (See
for example \cite[Exercise 6.19]{RW98}.) Let $w\in N_{N_{D}(\bar{x})}(-\nabla f(\bar{x}))\backslash\{0\}$.
Firstly, $-\lambda\nabla f(\bar{x})\in N_{D}(\bar{x})$ for all $\lambda\geq0$,
so $w^{T}(-\lambda\nabla f(\bar{x})-(-\nabla f(\bar{x})))\leq0$ for
all $\lambda\geq0$, that is $(1-\lambda)w^{T}\nabla f(\bar{x})\leq0$
for all $\lambda\geq0$. This implies that $w^{T}\nabla f(\bar{x})=0$,
or $w\in\{\nabla f(\bar{x})\}^{\perp}$. Secondly, $w\in N_{N_{D}(\bar{x})}(-\nabla f(\bar{x}))$
means that $w^{T}(s+\nabla f(\bar{x}))\leq0$ for all $s\in N_{D}(\bar{x})$.
Since $w^{T}\nabla f(\bar{x})=0$, this means that $w^{T}s\leq0$
for all $s\in N_{D}(\bar{x})$, which means that $w\in T_{D}(\bar{x})$.
Therefore, $w\in[\{\nabla f(\bar{x})\}^{\perp}\cap T_{D}(\bar{x})]\backslash\{0\}$,
which gives us the equivalence between the two conditions.
\end{proof}
In Proposition \ref{pro:Condn-Phi} below, we show that the conditions
in Lemma \ref{lem:strict-loc-min} can give us boundedness information
on the robust problem. We also give conditions for which $\bar{\Gamma}$
is a singleton.
\begin{prop}
\label{pro:Condn-Phi}(Conditions for $\bar{\Gamma}$ and optimality)
Suppose that $Q$ is Clarke regular at $\bar{A}\bar{x}-\bar{b}$,
and there is a vector $\gamma^{\prime}$ such that $\bar{A}\gamma^{\prime}\in\intr(T_{Q}(\bar{A}\bar{x}-\bar{b}))$.
Let $f$ be locally Lipschitz and subdifferentially regular at $\bar{x}$,
and $\bar{x}$ be a minimizer of the tangential problem \eqref{eq:FO_approx_1}.
Recall also $\Phi$, $\Psi$, $\Gamma$ and $\bar{\Gamma}$ as defined
in \eqref{eq:Phi}, \eqref{eq:Gamma} and \eqref{eq:barGamma}.
\begin{enumerate}
\item [(a)]The set $\bar{\Gamma}$ is bounded if there is some $c\in\partial f(\bar{x})$
such that $-c\in\intr(N_{\Phi}(\bar{x}))$. 
\item [(b)]If $f$ is $\mathcal{C}^{1}$ at $\bar{x}$ and $\nabla f(\bar{x})\neq0$,
then $\bar{\Gamma}$ is bounded if and only if $-\nabla f(\bar{x})\in\intr(N_{\Phi}(\bar{x}))$,
which is also equivalent to $\{\nabla f(\bar{x})\}^{\perp}\cap T_{\Phi}(\bar{x})=\{0\}$.
\item [(c)]A feasible $\bar{\gamma}\in\Gamma$ is in $\bar{\Gamma}$ if
we can find $c\in\partial f(\bar{x})$ such that $-c\in N_{T_{\Phi}(\bar{x})}(\bar{\gamma})$.
The condition $-c\in N_{T_{\Phi}(\bar{x})}(\bar{\gamma})$ holds when
we can find $\{u_{i}\}_{i=1}^{k}$ and $\{v_{i}\}_{i=1}^{k}$ such
that $v_{i}\in N_{T_{Q}(\bar{A}\bar{x}-\bar{b})}(\bar{\gamma}+u_{i})$,
$u_{i}\in L(\Delta\mathcal{U})$ and $-c=\sum_{i=1}^{k}\lambda_{i}\bar{A}^{T}v_{i}$
for some $\lambda_{i}\geq0$, $1\leq i\leq k$. Here, $L$ is the
linear map $L(\Delta A,\Delta b)=(\Delta A)\bar{x}-\Delta b$.
\item [(d)]The set $\bar{\Gamma}$ is a singleton if we can find some
$c\in\partial f(\bar{x})$ such that $-c\in\intr(N_{T_{\Phi}(\bar{x})}(\bar{\gamma}))$
for some $\bar{\gamma}\in\bar{\Gamma}$. The condition $-c\in\intr(N_{T_{\Phi}(\bar{x})}(\bar{\gamma}))$
holds when we can find $\{u_{i}\}_{i=1}^{k}$ and $\{v_{i}\}_{i=1}^{k}$
such that $v_{i}\in N_{T_{Q}(\bar{A}\bar{x}-\bar{b})}(\bar{A}\bar{\gamma}+u_{i})$,
$u_{i}\in L(\Delta\mathcal{U})$, $k=\dim(Q)$, $\{\bar{A}^{T}v_{i}\}_{i=1}^{k}$
is linearly independent, and $-c=\sum_{i=1}^{k}\lambda_{i}\bar{A}^{T}v_{i}$,
where $\lambda_{i}>0$ for all $1\leq i\leq k$. 
\end{enumerate}
\end{prop}
\begin{proof}
\textbf{Part (a):}  Seeking a contradiction, suppose that $\Gamma$
is unbounded, so there is a sequence $\{\gamma_{i}\}$ of solutions
of minimizers of \eqref{eq:FO_approx_1} such that $|\gamma_{i}|\to\infty$,
with $\frac{\gamma_{i}}{|\gamma_{i}|}\to\gamma_{\infty}$. Note \begin{eqnarray*}
\gamma_{i} & \in & \Gamma\\
 & = & \{\gamma:\bar{A}\gamma+L(\Delta\mathcal{U})\subset T_{Q}(\bar{A}\bar{x}-\bar{b})\}\\
 & \subset & \{\gamma:\bar{A}\gamma\in T_{Q}(\bar{A}\bar{x}-\bar{b})\}\\
 & = & T_{\Phi}(\bar{x}),\end{eqnarray*}
so $\gamma_{\infty}\in T_{\Phi}(\bar{x})$. On the other hand, since
$df(\bar{x})(\cdot)$ is continuous and positively homogeneous, we
have \begin{eqnarray*}
df(\bar{x})(\gamma_{\infty}) & = & \lim_{i\to\infty}df(\bar{x})\Big(\frac{\gamma_{i}}{|\gamma_{i}|}\Big)\\
 & = & \lim_{i\to\infty}\frac{1}{|\gamma_{i}|}df(\bar{x})(\gamma_{i})\\
 & = & 0.\end{eqnarray*}
The well known characterization of the subderivative in terms of support
functions (see for example \cite{Cla83,R70,RW98}) gives \begin{eqnarray*}
\max_{v\in\partial f(\bar{x})}v^{T}\gamma_{\infty} & = & df(\bar{x})(\gamma_{\infty})=0\\
\Rightarrow c^{T}\gamma_{\infty} & \leq & 0.\end{eqnarray*}
But $-c\in\intr(N_{\Phi}(\bar{x}))$, so $c^{T}\gamma_{\infty}>0$.
This contradiction tells us that the set $\bar{\Gamma}$ is bounded
as needed. 

\textbf{Part (b):} In view of part (a) and Lemma \ref{lem:equiv-strict},
we just need to prove that if $\bar{\Gamma}$ is bounded, then $\{\nabla f(\bar{x})\}^{\perp}\cap T_{\Phi}(\bar{x})=\{0\}$.
We prove the contrapositive. Suppose $w\in[\{\nabla f(\bar{x})\}^{\perp}\cap T_{\Phi}(\bar{x})]\backslash\{0\}$.
Then $\nabla f(\bar{x})^{T}w=0$, and $w\in T_{\Phi}(\bar{x})$, that
is $\bar{A}w\in T_{Q}(\bar{A}\bar{x}-\bar{b})$ by Proposition \ref{pro:T_Phi}.
Let $\bar{\gamma}$ be some element in $\bar{\Gamma}$. Then \begin{eqnarray*}
\bar{A}(\lambda w+\bar{\gamma})+L(\Delta\mathcal{U}) & = & \bar{A}\lambda w+[\bar{A}\bar{\gamma}+L(\Delta\mathcal{U})]\\
 & \subset & T_{Q}(\bar{A}\bar{x}-\bar{b})+T_{Q}(\bar{A}\bar{x}-\bar{b})\\
 & = & T_{Q}(\bar{A}\bar{x}-\bar{b}).\end{eqnarray*}
Therefore $\bar{\gamma}+\lambda w\in\bar{\Gamma}$ for all $\lambda\geq0$,
which shows that $\bar{\Gamma}$ is unbounded, concluding our proof.

\textbf{Parts (c), (d): }Note that \begin{eqnarray*}
\Psi & = & \{\gamma:\gamma+u\in T_{Q}(\bar{A}\bar{x}-\bar{b})\mbox{ for all }u\in L(\Delta\mathcal{U})\}\\
 & \subset & \{\gamma:\gamma+u_{i}\in T_{Q}(\bar{A}\bar{x}-\bar{b})\mbox{ for }1\leq i\leq k\}\\
 & = & \bigcap_{i=1}^{k}[T_{Q}(\bar{A}\bar{x}-\bar{b})-u_{i}],\end{eqnarray*}
so \begin{eqnarray*}
N_{\Psi}(\bar{A}\bar{\gamma}) & \supset & N_{\cap_{i=1}^{k}[T_{Q}(\bar{A}\bar{x}-\bar{b})-u_{i}]}(\bar{A}\bar{\gamma})\\
 & \supset & [N_{T_{Q}(\bar{A}\bar{x}-\bar{b})-u_{1}}(\bar{A}\bar{\gamma})]+\cdots+[N_{T_{Q}(\bar{A}\bar{x}-\bar{b})-u_{k}}(\bar{A}\bar{\gamma})]\\
 &  & \qquad\mbox{ (by [\ref{RW98}, Theorem 6.42])}\\
 & = & [N_{T_{Q}(\bar{A}\bar{x}-\bar{b})}(\bar{A}\bar{\gamma}+u_{1})]+\cdots+[N_{T_{Q}(\bar{A}\bar{x}-\bar{b})}(\bar{A}\bar{\gamma}+u_{k})]\\
 & \supset & \left\{ \sum_{i=1}^{k}\lambda_{i}v_{i}\mid\lambda_{i}\geq0\right\} .\end{eqnarray*}
Therefore, by Proposition \ref{pro:normal-of-tangent-feasible}, \begin{eqnarray*}
N_{\Gamma}(\bar{\gamma}) & = & \{\bar{A}^{T}y:y\in N_{\Psi}(\bar{A}\bar{\gamma})\}\\
 & \supset & \left\{ \sum_{i=1}^{k}\lambda_{i}\bar{A}^{T}v_{i}:\lambda_{i}\geq0\right\} .\end{eqnarray*}
For part (c), the condition stated is equivalent to the existence
of $c\in\partial f(\bar{x})$ such that $-c\in\{\sum_{i=1}^{k}\lambda_{i}v_{i}\mid\lambda_{i}\geq0\}$,
which implies $-c\in N_{\Gamma}(\bar{\gamma})$, which in turn implies
$\bar{\gamma}\in\bar{\Gamma}$. Part (d) follows by applying Lemma
\ref{lem:strict-loc-min}. 
\end{proof}
The conditions (c) and (d) in Proposition \ref{pro:Condn-Phi} can
be helpful for designing numerical methods for solving the tangential
problem. Due to Clarke regularity, the tangential problem is convex.
However, the problem of determining whether a point is feasible is
not necessarily easy.

The result corresponding to conditions (c) and (d) in Proposition
\ref{pro:Condn-Phi} can also be generalized for robust optimization
in general. The following result on normal cones in robust optimization
combined with results on optimality of nonlinear programs (in Lemma
\ref{lem:strict-loc-min} for example) gives us the optimality conditions.
The proof of the following result is a direct application of \cite[Theorems 6.14 and 6.42]{RW98}
similar to the proofs of Propositions \ref{pro:normal-of-tangent-feasible}
and \ref{pro:Condn-Phi}, so we shall only state the result.
\begin{prop}
(Normal cones in robust optimization) For $f:\mathbb{R}^{n}\to\mathbb{R}$,
consider the robust optimization problem\begin{eqnarray*}
 &  & \min_{x}f(x)\\
 &  & \mbox{s.t. }[\bar{A}+\Delta A]x-[\bar{b}-\Delta b]\in Q\mbox{ for all }(\Delta A,\Delta b)\in\Delta\mathcal{U}.\end{eqnarray*}
Let the sets $\Omega(\Delta A,\Delta b)\subset\mathbb{R}^{n}$ and
$\Xi\subset\mathbb{R}^{n}$ be defined by \begin{eqnarray*}
\Omega(\Delta A,\Delta b) & := & \{x\mid[\bar{A}+\Delta A]x-[\bar{b}-\Delta b]\in Q\},\\
\mbox{ and }\Xi & := & \{x\mid[\bar{A}+\Delta A]x-[\bar{b}-\Delta b]\in Q\mbox{ for all }(\Delta A,\Delta b)\in\Delta\mathcal{U}\}\\
 & = & \bigcap_{(\Delta A,\Delta b)\in\Delta\mathcal{U}}\Omega(\Delta A,\Delta b).\end{eqnarray*}
Let $\bar{x}\in\Xi$.
\begin{enumerate}
\item If $Q$ is Clarke regular at $[\bar{A}+\Delta A]\bar{x}-[\bar{b}-\Delta b]$
and the only vector $y\in N_{Q}([\bar{A}+\Delta A]\bar{x}-[\bar{b}-\Delta b])$
for which $[\bar{A}+\Delta A]^{T}y=0$ is $y=0$, then \[
N_{\Omega(\Delta A,\Delta b)}(\bar{x})=\{[\bar{A}+\Delta A]^{T}y\mid y\in N_{Q}([\bar{A}+\Delta A]\bar{x}-[\bar{b}-\Delta b])\}.\]

\item For any finite set $\{(\Delta A_{i},\Delta b_{i})\}_{i\in I}\subset\Delta\mathcal{U}$
such that $\bar{x}$ is Clarke regular at $\Omega(\Delta A_{i},\Delta b_{i})$
for all $i\in I$, the normal cone $N_{\Xi}(\bar{x})$ satisfies $N_{\Xi}(\bar{x})\supset\sum_{i\in I}N_{\Omega(\Delta A_{i},\Delta b_{i})}(\bar{x})$. 
\end{enumerate}
\end{prop}

\section{\label{sec:uncertain-sets}Addition of uncertainty sets in the tangential
problem}

For much of this section we focus on the tangential robust problem
on addition of uncertainty sets. More specifically, we ask what we
can say about the tangential problem with uncertainty set $\lambda_{1}\Delta\mathcal{S}_{1}+\lambda_{2}\Delta\mathcal{S}_{2}$
given knowledge of the optimal solutions of the tangential problem
for the uncertainty sets $\Delta\mathcal{S}_{1}$ and $\Delta\mathcal{S}_{2}$.
Such a problem can arise from having to considering robust optimization
problems with errors which are a sum of two or more unknown sources.

We begin with some elementary properties.
\begin{prop}
\label{pro:Trans+incl}(Elementary properties of uncertainty sets)
Suppose $Q$ is Clarke regular at $\bar{A}\bar{x}-\bar{b}$. For an
uncertainty set $\Delta\mathcal{U}$, suppose the solution of the
tangential problem is defined by \begin{eqnarray*}
v(\Delta\mathcal{U}) & := & \min_{\gamma}h(\gamma)\\
 &  & \mbox{s.t. }\bar{A}\gamma+[(\Delta A)\bar{x}-\Delta b]\in T_{Q}(\bar{A}\bar{x}-\bar{b})\mbox{ for all }(\Delta A,\Delta b)\in\Delta\mathcal{U},\end{eqnarray*}
where $h$ is sublinear. Then
\begin{itemize}
\item If $\Delta\mathcal{S}^{\prime}=\lambda\Delta\mathcal{S}$, then $v(\Delta\mathcal{S}^{\prime})=\lambda v(\Delta\mathcal{S})$.
\item If $\Delta\mathcal{S}^{\prime}\supset\Delta\mathcal{S}$, then $v(\Delta\mathcal{S})\leq v(\Delta\mathcal{S}^{\prime})$.
\item If $h(\gamma)=c^{T}\gamma$ and $\Delta\mathcal{S}^{\prime}=\Delta\mathcal{S}+\{\bar{s}\}$,
then $v(\Delta\mathcal{S}^{\prime})=v(\Delta\mathcal{S})+c^{T}\bar{s}$.
\end{itemize}
\end{prop}
We have the following result to study how set addition affects the
solution to the tangential problem.
\begin{prop}
\label{pro:set-add}(Set addition) Recall the definition of $v(\Delta\mathcal{U})$
in Proposition \ref{pro:Trans+incl}. Suppose $\Delta\mathcal{S}=\lambda_{1}\Delta\mathcal{S}_{1}+\lambda_{2}\Delta\mathcal{S}_{2}$,
where $\lambda_{1},\lambda_{2}\geq0$. Then $v(\Delta\mathcal{S})\leq\lambda_{1}v(\Delta\mathcal{S}_{1})+\lambda_{2}v(\Delta\mathcal{S}_{2})$. \end{prop}
\begin{proof}
In view of Proposition \ref{pro:Trans+incl}, we only need to prove
the case for $\lambda_{1}=\lambda_{2}=1$. Suppose $\bar{\gamma}_{1}$
and $\bar{\gamma}_{2}$ are solutions to \begin{eqnarray*}
 &  & \min_{\gamma}h(\gamma)\\
 &  & \mbox{s.t. }\bar{A}\gamma+[(\Delta A)\bar{x}-\Delta b]\in T_{Q}(\bar{A}\bar{x}-\bar{b})\mbox{ for all }(\Delta A,\Delta b)\in\Delta\mathcal{S}_{i}\end{eqnarray*}
for $i=1,2$. If $(\Delta A,\Delta b)\in\Delta\mathcal{S}_{1}+\Delta\mathcal{S}_{2}$,
then \[
[(\Delta A)\bar{x}-\Delta b]=[(\Delta A_{1})\bar{x}-\Delta b_{1}]+[(\Delta A_{1})\bar{x}-\Delta b_{2}],\]
for some $[(\Delta A_{i})\bar{x}-\Delta b_{i}]\in\Delta\mathcal{S}_{i}$
for $i=1,2$. This means that \begin{eqnarray*}
 &  & \bar{A}(\bar{\gamma}_{1}+\bar{\gamma}_{2})-[(\Delta A)\bar{x}-\Delta b]\\
 & = & \big[\bar{A}\bar{\gamma}_{1}-[(\Delta A_{1})\bar{x}-\Delta b_{1}]\big]+\big[\bar{A}\bar{\gamma}_{2}-[(\Delta A_{2})\bar{x}-\Delta b_{2}]\big]\\
 & \in & T_{Q}(\bar{A}\bar{x}-\bar{b})+T_{Q}(\bar{A}\bar{x}-\bar{b})\\
 & = & T_{Q}(\bar{A}\bar{x}-\bar{b}).\end{eqnarray*}
So $\bar{\gamma}_{1}+\bar{\gamma}_{2}$ is a feasible, though not
necessarily optimal, solution to the tangential problem where the
uncertainty set is $\Delta\mathcal{S}_{1}+\Delta\mathcal{S}_{2}$,
which shows that $v(\Delta\mathcal{S})\leq h(\bar{\gamma}_{1}+\bar{\gamma}_{2})\leq h(\bar{\gamma}_{1})+h(\bar{\gamma}_{2})=v(\Delta\mathcal{S}_{1})+v(\Delta\mathcal{S}_{2})$. 
\end{proof}
In a nondegenerate linear programming problem, we do have equality
in Proposition \ref{pro:set-add}. The assumption that $\bar{A}$
is an invertible square matrix below is not restrictive, because this
is exactly what happens in a nondegenerate linear program.
\begin{prop}
\label{pro:set-add-nice}(Set addition in nondegenerate linear programming)
Consider the tangential problem having a linear programming structure\begin{eqnarray*}
v(\Delta\mathcal{U}) & := & \min_{\gamma}c^{T}\gamma\\
 &  & \mbox{s.t. }\bar{A}\gamma+[(\Delta A)\bar{x}-\Delta b]\leq0\mbox{ for all }(\Delta A,\Delta b)\in\Delta\mathcal{U}.\end{eqnarray*}
Assume that $\bar{A}$ is square and invertible, and $-c=\sum_{i=1}^{k}\lambda_{i}\bar{A}_{i}$
for some $\lambda_{i}>0$, $1\leq i\leq k$. Then $v(\Delta\mathcal{S}_{1}+\Delta\mathcal{S}_{2})=v(\Delta\mathcal{S}_{1})+v(\Delta\mathcal{S}_{2})$.\end{prop}
\begin{proof}
Recall that feasibility can be rewritten as \[
\bar{A}_{i}\gamma+\max_{(\Delta A,\Delta b)\in\Delta\mathcal{U}}[(\Delta A_{i})\bar{x}-\Delta b_{i}]\leq0\mbox{ for all }i.\]
Form the vector $w(\Delta\mathcal{U})\in\mathbb{R}^{m}$ by \[
[w(\Delta\mathcal{U})]_{i}:=-\max_{(\Delta A,\Delta b)\in\Delta\mathcal{U}}[(\Delta A_{i})\bar{x}-\Delta b_{i}].\]
The value $v(\Delta\mathcal{U})$ is equal to $c^{T}\bar{A}^{-1}w(\Delta\mathcal{U})$.
The condition on $-c$ in the statement assures that this minimizer
is unique through convexity and Lemma \ref{lem:strict-loc-min}. It
is clear that for each $i$, we have \begin{eqnarray*}
 &  & \max_{(\Delta A,\Delta b)\in\Delta\mathcal{S}_{1}+\Delta\mathcal{S}_{2}}[(\Delta A_{i})\bar{x}-\Delta b_{i}]\\
 & = & \max_{(\Delta A,\Delta b)\in\Delta\mathcal{S}_{1}}[(\Delta A_{i})\bar{x}-\Delta b_{i}]+\max_{(\Delta A,\Delta b)\in\Delta\mathcal{S}_{2}}[(\Delta A_{i})\bar{x}-\Delta b_{i}],\end{eqnarray*}
which implies that $w(\Delta\mathcal{S}_{1}+\Delta\mathcal{S}_{2})=w(\Delta\mathcal{S}_{1})+w(\Delta\mathcal{S}_{2})$.
This immediately gives $v(\Delta\mathcal{S}_{1}+\Delta\mathcal{S}_{2})=v(\Delta\mathcal{S}_{1})+v(\Delta\mathcal{S}_{2})$
as needed. 
\end{proof}
The next step is to ask whether the property in Proposition \ref{pro:set-add-nice}
is satisfied for problems of the form\begin{eqnarray*}
 &  & \min_{\gamma}c^{T}\gamma\\
 &  & \mbox{s.t. }\bar{A}\gamma+[(\Delta A)\bar{x}-\Delta b]\in Q\mbox{ for all }(\Delta A,\Delta b)\in\Delta\mathcal{U},\end{eqnarray*}
where $Q$ is a convex cone. If $Q=\{y\mid\tilde{A}y\leq0\}$, where
$\tilde{A}$ is an invertible square matrix, then the constraint above
can be transformed into \[
\tilde{A}\big[\bar{A}\gamma+[(\Delta A)\bar{x}-\Delta b]\big]\leq0\mbox{ for all }(\Delta A,\Delta b)\in\Delta\mathcal{U},\]
and we can apply Proposition \ref{pro:set-add-nice} on $\tilde{A}\bar{A}$. 

\begin{figure}[h]
\begin{tabular}{|c|c|}
\hline 
\includegraphics[scale=0.4]{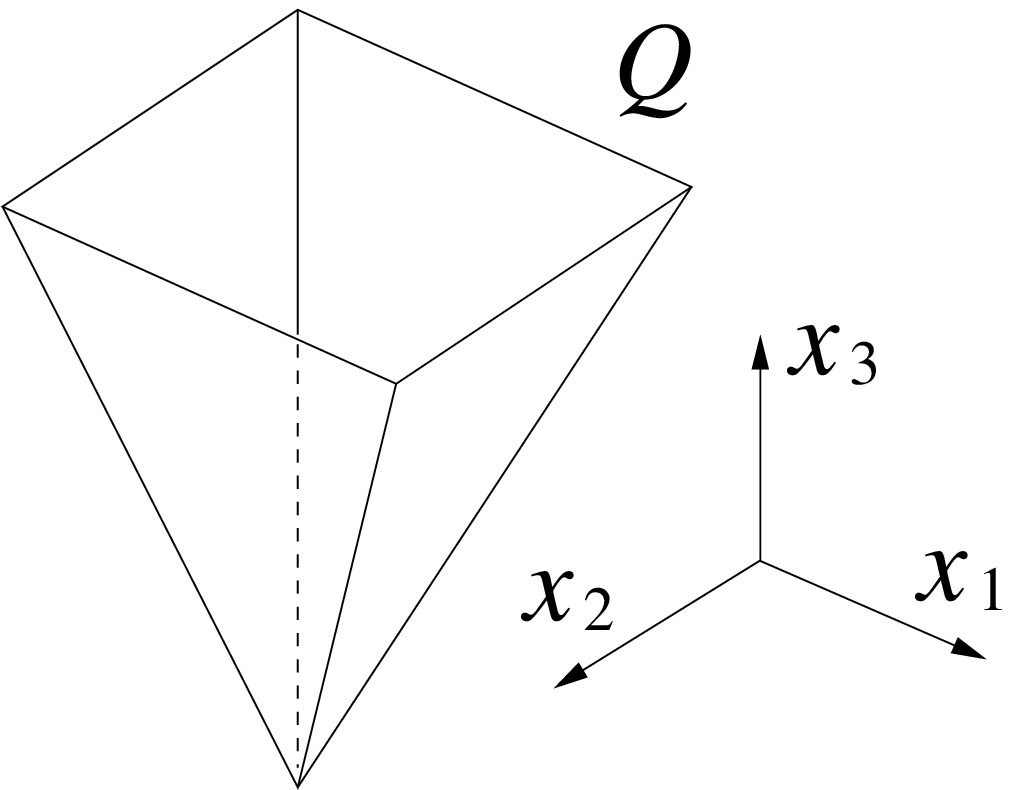} & \includegraphics[scale=0.4]{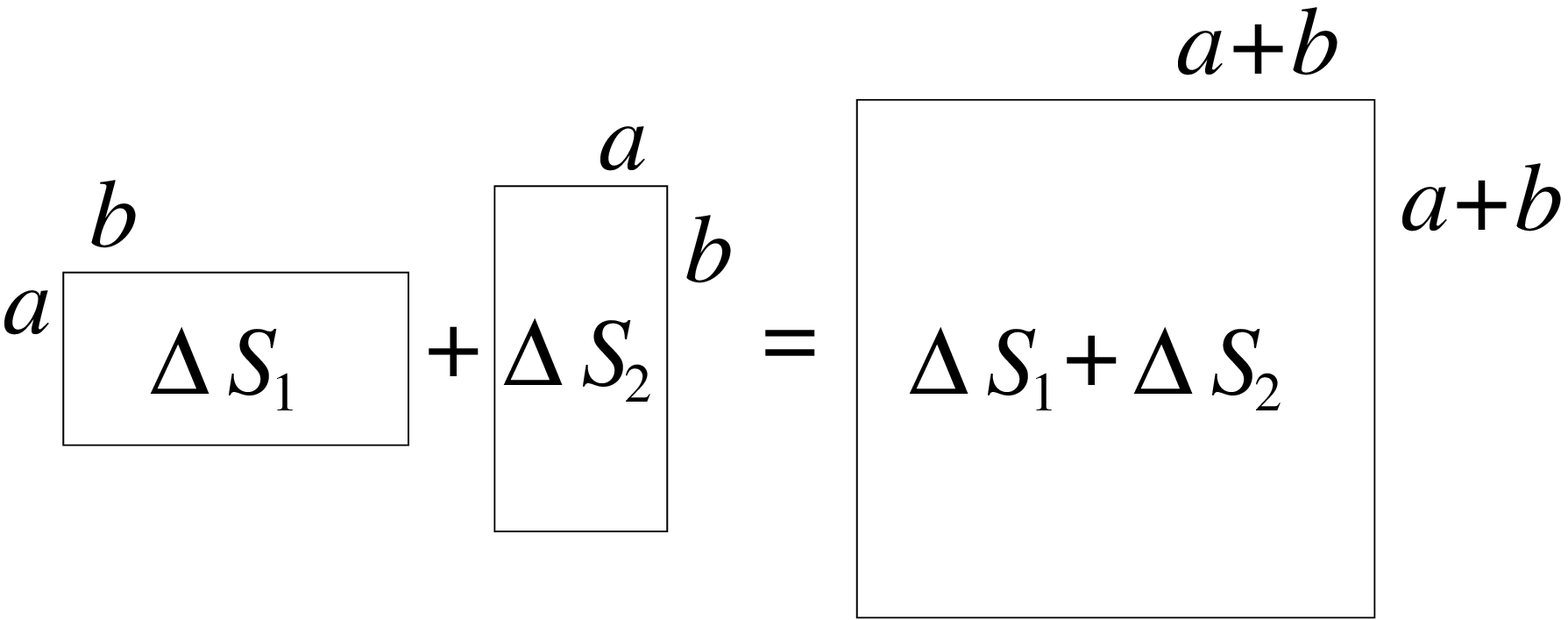}\tabularnewline
\hline
\end{tabular}

\caption{\label{fig:square-eg}These diagrams illustrate Example \ref{exa:square-eg}.
The figure on the left illustrates the convex cone $Q$ in $\mathbb{R}^{3}$,
while the figure on the right illustrates the first 2 coordinates
for the shapes $\Delta\mathcal{S}_{1}$ and $\Delta\mathcal{S}_{2}$.}

\end{figure}

In the general case, the equality $v(\Delta\mathcal{S}_{1}+\Delta\mathcal{S}_{2})=v(\Delta\mathcal{S}_{1})+v(\Delta\mathcal{S}_{2})$
may not hold. We give two examples to illustrate this. In the first
example, we have a degenerate linear program, while in the second
example, we have different cones for a conic programming problem.
\begin{example}
\label{exa:degen}(Inequality in sums of uncertainty sets 1) Consider
the problem \begin{eqnarray*}
\bar{v}(\Delta\mathcal{U}) & := & \min_{x\in\mathbb{R}^{2}}x_{2}\\
 &  & \mbox{s.t. }(\bar{A}+\Delta A)x-(\mathbf{0}+\Delta b)\leq0\mbox{ for all }(\Delta A,\Delta b)\in\Delta\mathcal{U},\\
 &  & \mbox{,where }\bar{A}=\left(\begin{array}{cc}
-1 & -1\\
1 & -1\\
0.5 & -1\end{array}\right)\mbox{, }\bar{b}=\mathbf{0}.\end{eqnarray*}
The optimal solution to the nonrobust problem is $\bar{x}=0$, and
the tangential problem is\begin{eqnarray*}
v(\Delta\mathcal{U}) & := & \min_{\gamma\in\mathbb{R}^{2}}\gamma_{2}\\
 &  & \mbox{s.t. }\bar{A}\gamma+[(\Delta A)\mathbf{0}-\Delta b]\leq0\mbox{ for all }(\Delta A,\Delta b)\in\Delta\mathcal{U}.\end{eqnarray*}
We illustrate this example in Figure \ref{fig:degen}. Let the uncertainty
sets $\Delta\mathcal{S}_{1}$ and $\Delta\mathcal{S}_{2}$ be defined
by \begin{eqnarray*}
\Delta\mathcal{S}_{1} & := & \{\mathbf{0}\}\times\{\Delta b:\Delta b_{1}=0,|\Delta b_{2}|\leq4\epsilon,\Delta b_{3}=0\}\\
\Delta\mathcal{S}_{2} & := & \{\mathbf{0}\}\times\{\Delta b:\Delta b_{1}=0,\Delta b_{2}=0,|\Delta b_{3}|\leq3\epsilon\}\\
\Delta\mathcal{S}_{1}+\Delta\mathcal{S}_{2} & = & \{\mathbf{0}\}\times\{\Delta b:\Delta b_{1}=0,|\Delta b_{2}|\leq4\epsilon,|\Delta b_{3}|\leq3\epsilon\}.\end{eqnarray*}
To find $v(\Delta\mathcal{S}_{1})$, we note that tangential problem
becomes \[
\left(\begin{array}{cc}
-1 & -1\\
1 & -1\\
0.5 & -1\end{array}\right)\gamma\leq\left(\begin{array}{c}
0\\
-4\epsilon\\
0\end{array}\right).\]
The first two rows are the active constraints, which gives a solution
of $\bar{\gamma}(\Delta\mathcal{S}_{1})=(-2\epsilon,2\epsilon)$ and
$v(\Delta\mathcal{S}_{1})=2\epsilon$. To find $v(\Delta\mathcal{S}_{2})$,
a similar set of calculations shows that the first and third constraints
are the active constraints, which gives $\bar{\gamma}(\Delta\mathcal{S}_{2})=(-2\epsilon,2\epsilon)$
and $v(\Delta\mathcal{S}_{2})=2\epsilon$. Similarly, $v(\Delta\mathcal{S}_{1}+\Delta\mathcal{S}_{2})=2\epsilon$,
and all constraints are active. We have $v(\Delta\mathcal{S}_{1}+\Delta\mathcal{S}_{2})<v(\Delta\mathcal{S}_{1})+v(\Delta\mathcal{S}_{2})$
as needed. 

\begin{figure}
\includegraphics[scale=0.5]{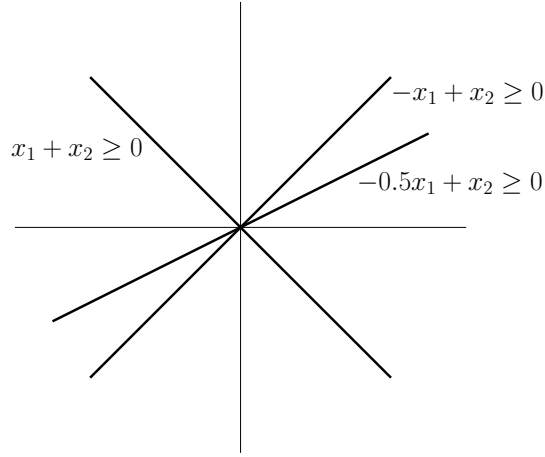}

\caption{\label{fig:degen}Illustration of Example \ref{exa:degen}. }

\end{figure}

\end{example}
Here is a second example for the case when the cone $Q$ is slightly
more complicated.
\begin{example}
\label{exa:square-eg}(Inequality in sums of uncertainty sets 2) Consider
the problem \begin{eqnarray*}
\bar{v}(\Delta\mathcal{U}) & := & \min_{\gamma\in\mathbb{R}^{3}}x_{3}\\
 &  & \mbox{s.t.}(I+\Delta A)x-(\mathbf{0}+\Delta b)\in Q\mbox{ for all }(\Delta A,\Delta b)\in\Delta\mathcal{U}.\end{eqnarray*}
Here, the matrix $\bar{A}=I$ is the $3\times3$ identity matrix,
$\bar{b}\in\mathbb{R}^{3}$ is the zero vector, and the convex cone
$Q$ is defined by\[
Q:=\{x\in\mathbb{R}^{3}:x_{3}\geq\max(|x_{1}|,|x_{2}|)\}.\]
The optimal solution to the nonrobust problem is $\bar{x}=\mathbf{0}$,
and the tangential problem is \begin{eqnarray*}
v(\Delta\mathcal{U}) & := & \min_{\gamma\in\mathbb{R}^{3}}\gamma_{3}\\
 &  & \mbox{s.t.}\gamma+[(\Delta A)\mathbf{0}-\Delta b]\in T_{Q}(\mathbf{0})\mbox{ for all }(\Delta A,\Delta b)\in\Delta\mathcal{U}.\end{eqnarray*}
Let $a<b$ and the sets $\mathcal{S}_{1}$ and $\mathcal{S}_{2}$
be defined by \begin{eqnarray*}
\Delta\mathcal{S}_{1} & := & \{\mathbf{0}\}\times\Big\{\Delta b\in\mathbb{R}^{3}:|\Delta b_{1}|\leq\frac{a}{2},|\Delta b_{2}|\leq\frac{b}{2},|\Delta b_{3}|\leq\frac{\delta}{2}\Big\}\\
\Delta\mathcal{S}_{2} & := & \{\mathbf{0}\}\times\Big\{\Delta b\in\mathbb{R}^{3}:|\Delta b_{1}|\leq\frac{b}{2},|\Delta b_{2}|\leq\frac{a}{2},|\Delta b_{3}|\leq\frac{\delta}{2}\Big\}\\
\Delta\mathcal{S}_{1}+\Delta\mathcal{S}_{2} & = & \{\mathbf{0}\}\times\Big\{\Delta b\in\mathbb{R}^{3}:|\Delta b_{1}|\leq\frac{a+b}{2},|\Delta b_{2}|\leq\frac{a+b}{2},|\Delta b_{3}|\leq\delta\Big\}.\end{eqnarray*}
See Figure \ref{fig:square-eg} for an illustration of the convex
cone $Q$, and the projection of $\Delta\mathcal{S}_{1}$, $\Delta\mathcal{S}_{2}$
and $\Delta\mathcal{S}_{1}+\Delta\mathcal{S}_{2}$ onto the 2-dimensional
space corresponding to the first 2 coordinates in $\mathbb{R}^{3}$.
It is elementary that $v(\Delta\mathcal{S}_{1})=v(\Delta\mathcal{S}_{2})=\frac{b}{2}+\frac{\delta}{2}$,
and $v(\Delta\mathcal{S}_{1}+\Delta\mathcal{S}_{2})=\frac{a+b}{2}+\delta$.
This gives $v(\Delta\mathcal{S}_{1}+\Delta\mathcal{S}_{2})<v(\Delta\mathcal{S}_{1})+v(\Delta\mathcal{S}_{2})$.
If the calculations had been performed with the second order cone
$Q^{\prime}$ defined by \[
Q^{\prime}:=\Big\{x\in\mathbb{R}^{3}:x_{3}\geq\sqrt{x_{1}^{2}+x_{2}^{2}}\Big\}\]
instead, then we also get the conclusion $v(\Delta\mathcal{S}_{1}+\Delta\mathcal{S}_{2})<v(\Delta\mathcal{S}_{1})+v(\Delta\mathcal{S}_{2})$.
\end{example}

\section*{Acknowledgments}

I am grateful to Henry Wolkowicz for conversations and his initial
ideas in \cite{WY09} that led to me studying this problem, and to
the Natural Sciences Engineering Research Council (Canada) for supporting
the research. The research for this paper was carried out in the Department
of Combinatorics and Optimization, University of Waterloo, when I
was on a postdoctoral position there, and I gratefully acknowledge
them for providing a splendid working environment.


\begin{thebibliography}{10}
\bibitem{AG03}F. Alizadeh and D. Goldfarb, \emph{Second-order cone
programming,} Mathematical Programming, Ser B 95:3--51 (2003).

\bibitem{AW81}H. Attouch and R. J.-B. Wets, \emph{Approximation and
convergence in nonlinear optimization}, in Nonlinear Programming 4,
edited by O. Mangasarian, R. Meyer and S. Robinson, pp. 367--395,
Academic Press, New York, 1981.

\bibitem{BGN09}A. Ben-Tal, L. El Ghaoui and A. Nemirovski. \emph{Robust
Optimization}, Princeton Series in Applied Mathematics, Princeton,
2009.

\bibitem{Cla83}F.H. Clarke. \emph{Optimization and Nonsmooth Analysis}.
Wiley, New York, 1983. Republished as Vol. 5, Classics in Applied
Mathematics, SIAM, 1990.

\bibitem{Mor06}B.S. Mordukhovich. \emph{Variational Analysis and
Generalized Differentiation I and II}., Grundlehren der mathematischen
Wissenschaften, Vols 330 \& 331, Springer, Berlin, 2006.

\bibitem{R70}R. T. Rockafellar, \emph{Convex Analysis}, Princeton,
1970.

\bibitem{R79}R. T. Rockafellar, \emph{Clarke's tangent cones and
the boundaries of closed sets in $\mathbb{R}^{n}$}, Nonlinear Analysis:
Theory, Methods and Applications, 3, 145--154, 1979.

\bibitem{RW84}R. T. Rockafellar and R. J.-B. Wets, \emph{Variational
Systems, an introduction}, in Multifunctions and Integrands: Stochastic
Analysis, Approximation and Optimization, edited by G. Salinetti,
Lecture Notes in Mathematics, 1091, pp. 1--54, Springer Verlag, Berlin,
1984.

\bibitem{RW98} \label{RW98}R.T. Rockafellar and R.J.-B. Wets. \emph{Variational
Analysis}, Grundlehren der mathematischen Wissenschaften, Vol 317,
Springer, Berlin, 1998.

\bibitem{WY09}H. Wolkowicz and W.L.N. Yeung,\emph{ An alternative
approach to sensitivity analysis in linear programming}, 2009 (unpublished).
\end{thebibliography}
\end{document}